\documentclass[a4paper,12pt]{amsart}

\usepackage{amssymb,amsmath,latexsym,color,graphicx}
\usepackage{hyperref}
\usepackage{amsfonts}
\usepackage{amsthm}
\usepackage[color]{showkeys}
\usepackage{tikz}
\usetikzlibrary{positioning, arrows,decorations.pathmorphing}
\begin{document}

\title[]{Multi-domain spectral approach with Sommerfeld condition for 
the Maxwell equations}
\author{Christian Klein}
\address{Institut de Math\'ematiques de Bourgogne, UMR 5584\\
                Universit\'e de Bourgogne-Franche-Comt\'e, 9 avenue Alain Savary, 21078 Dijon
                Cedex, France\\
    E-mail Christian.Klein@u-bourgogne.fr}

\author{Nikola Stoilov}
\address{Institut de Math\'ematiques de Bourgogne, UMR 5584\\
                Universit\'e de Bourgogne-Franche-Comt\'e, 9 avenue Alain Savary, 21078 Dijon
                Cedex, France\\
    E-mail Nikola.Stoilov@u-bourgogne.fr}
\date{\today}
\begin{abstract}
	We present a  multi-domain spectral approach with an exterior compactified 
	domain for the Maxwell equations for monochromatic fields. The Sommerfeld radiation condition is imposed exactly 
	at infinity being a finite point on the numerical grid. As  an  
	example, axisymmetric situations in spherical and prolate 
	spheroidal coordinates are discussed, as well 
	as the interaction of a radiating dipole with a nano-particle.  
	
\end{abstract}

\thanks{We thank the anonymous referees for helpful suggestions and 
remarks. This work is partially supported by 
the ANR-FWF project ANuI - ANR-17-CE40-0035, the isite BFC project 
NAANoD, the EIPHI Graduate School (contract ANR-17-EURE-0002), by the 
European Union Horizon 2020 research and innovation program under the 
Marie Sklodowska-Curie RISE 2017 grant agreement no. 778010 IPaDEGAN 
and the EITAG project funded by the FEDER de Bourgogne, the region 
Bourgogne-Franche-Comt\'e and the EUR EIPHI.}

\maketitle

\section{Introduction}
The interaction between electromagnetic radiation and matter is  
arguably one of the most important problems of physics, and one of great practical importance. The governing 
equations for this are the Maxwell equations, and their efficient 
numerical solution in situations appearing in applications is thus 
crucial. Interestingly the classical Maxwell equations are 
also relevant in the context of quantum emitters since the latter can be 
treated as a quantum system interacting with a classical field, see 
for instance \cite{DRJCCG,SHML} and references therein. In general 
one is not interested in the solution of a particular initial value 
problem in this context, but a discussion in the frequency domain. In 
this case, the Maxwell equations in non-magnetizable matter can be 
cast into the form (see section \ref{sec:max} for a 
short derivation and for references) of a vector Helmholtz equation,  
\begin{equation}
\nabla\times\nabla\times \mathbf{E}(\mathbf{x},\omega) - 
\omega^2\epsilon(\mathbf{r},\omega)\mathbf{E}(\mathbf{x},\omega) =  \mathbf{f}. 
\label{Helmholtz}
\end{equation} 
Here $\mathbf{x}\in \mathbb{R}^{3}$ with components $x_{i}$, 
$i=1,2,3$, $\mathbf{E}\in \mathbb{C}^{3}$ is the electric field, 
$\omega\in\mathbb{R}$ is 
the frequency, $\epsilon(\mathbf{r},\omega)$ is the 
permittivity in the Maxwell equations with matter, see 
(\ref{maxwell}) and (\ref{matter}), $\mathbf{f}$ is an inhomogenity due to free currents, 
and $\nabla$ is the vector operator with components 
$\nabla_{i}=\frac{\partial}{\partial x_{i}}$, $i=1,2,3$. It is the 
goal of this paper to provide a multi-domain spectral method for the 
solution of this equation for systems corresponding to a localized 
isolated matter configuration thus extending the method by Gharti et 
al \cite{GJ17}- \cite{GLT3}. 

In order to get a unique solution to the Helmholtz equation 
(\ref{Helmholtz}), a condition needs to be imposed at infinity. 
Sommerfeld suggested that there should be no incoming radiation at 
infinity, i.e., the only source of radiation should be the matter 
distribution. 
The Sommerfeld radiation condition \cite{SOM49} reads
\begin{equation}
	\lim_{||\mathbf{x}||\to\infty}||\mathbf{x}||\left(e^{i\omega||\mathbf{x}||}\mathbf{E}(\mathbf{x},\omega)-1\right) = 0
	\label{som},
\end{equation}
where $||\cdot||$ is the euclidean norm (in vacuum 
$\epsilon(\mathbf{x},\omega)=1$). This implies, however, that the 
solution has an oscillatory singularity (proportional 
to $e^{i\omega||\mathbf{x}||}$) for large $||\mathbf{x}||$ which is 
numerically challenging. What makes things worse is that it is known, 
see \cite{ATK,WIL}, that the solution of the scalar Helmholtz equation 
$$\Delta E + \omega^{2}E=0,$$ 
where $\Delta$ is the three-dimensional 
Laplace operator, with a 
Sommerfeld condition has in spherical coordinates (see section 
\ref{spher}) for large $r$ the 
form
\begin{equation}
	E = \frac{e^{-i\omega r}}{\omega 
	r}\sum_{n=0}^{\infty}\frac{a_{n}(\theta,\phi)}{(\omega r)^{n}},
	\label{asym}
\end{equation}
where the $a_{n}$, $n\in \mathbb{R}$ depend only on $\theta,\phi$. 
This means that the solution is not only  oscillatory near 
infinity, but also slowly decaying in $1/r$.  

Imposing boundary conditions at infinity has been discussed in many works and in various contexts. 
The most popular approach is to  truncate the problem and solve it on 
a finite domain by imposing artificial boundaries at a finite 
distance from the physical system (for a general review see e.g. 
\cite{GIV04, T98} and references therein).  One option in this case 
is to consider non-reflecting boundaries (NRBs), that is, boundaries 
that allow the waves to leave the truncated domain without spurious 
reflections that may pollute the solution in the computational domain 
of interest.  There are two main groups of NRBs, namely, 
Non-Reflecting Boundary Conditions (NRBCs) and Non-Reflecting 
Boundary Layers (NRBLs). NRBCs are boundary conditions on the 
artificial boundary that absorb incident waves, see for example 
\cite{BGT82} \cite{ECB14}. On the other hand, NRBLs are additional 
computational domains which absorb waves that are traveling inside 
the layer, effectuating trivial boundary condition at the end of the 
domain towards infinity. The most prominent among the NRBL techniques 
is the Perfectly Matched Layer (PML) initially  developed by 
B\'erenger in 1994 for electromagnetic scattering \cite{BER94}. The 
idea is to add an absorbing layer so that plane waves produce no 
reflection and that inside the layer the solution decays 
exponentially \cite{GIV08}. Such approaches work very well for linear 
problems as studied here, but need in general an optimization of the 
parameters of the absorbing layer, see for instance the discussion in 
\cite{birem}. On the other hand in the case of slowly decaying 
potentials as in electromagnetism and relativity  approximate implementation of 
non-reflecting boundary conditions can lead to non-negligible 
errors, see for instance the comparison of compactification (as in 
the present paper)
and truncation schemes in a relativistic context in \cite{eriguchi}. 
In a time dependent context, see \cite{SN} for an implementation of 
the Sommerfeld condition in this context, for compactification 
approaches based on the conformal invariance the review \cite{joerg}. 
Note that the techniques explained in this paper can be directly 
applied to similar problems in linearized gravity.

Another approach, which is closely connected to the one taken in the 
present work, is that of mapped infinite elements. Its origin can be 
traced back to the works of Zienkevich and Bettess \cite{ZKB79}, for 
a comprehensive review see \cite{BET92}. The basic idea in one 
dimension, $x$, is to add an element extending to infinity, where we 
map to a new coordinate so that $x= 2x_0/(1-\xi)$. The infinite 
element is  thus mapped onto $[-1,1]$ and infinity becomes the 
regular grid point $\xi = 1$.   In a general setting, shape functions  $M_i$ in the infinite element are polynomials in $(\xi)$, which translates to polynomials in $1/x$. The Sommerfeld condition can be directly implemented on the shape functions and take them in the form $M_i = e^{-i\omega x}P_i(\xi)$.  
The idea was further developed by Beer and Meeks in 
\cite{BM81}. It was applied under the name ``infinite boundary 
element''  for electromagnetic and other problems by Kagawa \emph{et 
al} \cite{KYK83}.  Recently this was  developed into `spectral infinite element' methods by Gharti {\emph et al.} \cite{GJ17, GTZ18} when dealing with geophysical problems. 

In this paper, we use an approach similar to the mapped infinite 
elements \cite{ZKB79,BET92} for of a multi-domain spectral 
approach. Note that spectral methods are distinguished by their 
excellent approximation properties for analytical functions since the 
numerical error  in such a case decreases exponentially with the 
numerical resolution. They are thus especially effective if the 
function $\epsilon(\mathbf{x},\omega)$ in (\ref{Helmholtz}) is 
analytical in the considered domains. We concentrate here on the case 
where this is true on concentric spheres or spheroids, for instance a 
spherical or spheroidal conductor in vacuum, possibly with multiple 
layers. More precisely, we consider a number 
$N_{d}+1$, $N_{d}\in 
\mathbb{N}$,  of domains where $\epsilon$ is smooth in each of them, 
but in concrete examples we discuss the case of three domains, an interior domain such that the boundary of the matter is a 
domain boundary, a second domain in vacuum\footnote{Experience in an 
astrophysical context in \cite{lorene} shows that it is numerically 
recommended not to match the infinite domain directly to the matter 
configuration, but to apply an intermediate domain.}, and a third domain 
with the local parameter $1/||\mathbf{x}||$ around infinity. As in 
\cite{ZKB79,BET92} we split off the oscillatory term in (\ref{asym}) 
by writing 
$$E = e^{-i\omega r}\tilde{E}$$
and solving the equations in the compactified domain for $\tilde{E}$ 
which is non-oscillatory and analytical in $1/||\mathbf{x}||$ and 
thus ideally suited for a  spectral method. 





The paper is organized as follows: in section \ref{sec:max}, we 
review the Maxwell equations in spherical and prolate spheroidal 
coordinates and introduce the twist potential in the axisymmetric case. 
In section \ref{sec:bound} we discuss the matching and 
the Sommerfeld radiation condition. In section 4 we present our
numerical approach. Some examples are discussed in section 5. We add 
some concluding remarks in section 6.

\textbf{Notation:} Partial derivatives of a function $u$ with 
respect to $x$ are denoted by $\partial_{x}u$ or $u_{x}$, vector 
indices are superscripts.

\section{The Maxwell equations in spherical and prolate spheroidal 
coordinates}\label{sec:max}
In this section we give a brief summary of the Maxwell equations in matter. 
A convenient form for a numerical 
solution is presented in the axisymmetric case in spherical and 
prolate spheroidal coordinates. 

Throughout this paper we assume that the studied problems allow 
a Fourier transform in $t$ with $\omega$ being the dual 
Fourier variable to $t$. This means we are interested in the long 
time interaction between electromagnetic radiation rather than in 
specific initial value problems. The Maxwell equations in this case 
read
\begin{align}
	\nabla\cdot \mathbf{D}(\mathbf{x},\omega)=\sigma ,\quad & \nabla\cdot \mathbf{H}(\mathbf{x},\omega)=0,
	\nonumber\\
	\nabla \times \mathbf{E}(\mathbf{x},\omega) = -i\omega \mathbf{B}(\mathbf{x},\omega), \quad & \nabla \times 
	\mathbf{H}(\mathbf{x},\omega) = \mathbf{J}+i\omega 
	\mathbf{D}(\mathbf{x},\omega),
	\label{maxwell}
\end{align}
where $\sigma$ is the density of the free charges, and $\mathbf{J}$ is the 
density of the free currents. Note that we use geometric units here 
in which the velocity of light is equal to 1. 
We assume that the matter is such that the following relations hold
\begin{equation}
	\mathbf{H}(\mathbf{x},\omega) = \mathbf{B}(\mathbf{x},\omega),\quad \mathbf{D}(\mathbf{x},\omega) = 
	\epsilon(\mathbf{x},\omega)\mathbf{E}(\mathbf{x},\omega).
	\label{matter}
\end{equation}
In this case the Maxwell equations (\ref{maxwell}) are equivalent to 
the vector Helmholtz equation (\ref{Helmholtz}). 

The permittivity $\epsilon(\mathbf{x},\omega)$ can depend on 
$\mathbf{x}$, and it is assumed here that it is piece-wise smooth. This 
means that it is a smooth function in various domains of 
$\mathbb{R}^{3}$, but not necessarily continuous on the whole of 
$\mathbb{R}^{3}$. A typical situation would be that it is constant in 
a homogeneous compact matter distribution, but may jump at the boundary of 
the matter which is assumed to be smooth. The matter has to be of finite extension which means that 
there is vacuum for $||\mathbf{x}||\gg 1$, i.e.,  
$\epsilon(\mathbf{x},\omega)=1$ for $||\mathbf{x}||\to\infty$. 

It is well known, see for instance \cite{KR} and references therein, 
that the Maxwell equations in the presence of a symmetry given by a 
Killing vector (in the axisymmetric case to be considered here this is $\partial_{\phi}$) simplify considerably. In this 
case one can introduce the \emph{twist potential}, here the $\phi$ 
component of $\nabla\times \mathbf{E}$. The Maxwell equations are 
then equivalent to a scalar equation for the twist potential. Whereas 
our numerical approach does not rely on the presence of such a 
symmetry, we limit ourselves to this case here since it allows for a 
more compact presentation of the approach. 

In applications, for instance in the 
description of nano-conductors in optics, see  
\cite{DRJCCG} and references therein, the Green's functions of the 
studied Maxwell equations play an important role. Therefore we 
mention how the latter are related to the equations we actually solve.

\subsection{Spherical coordinates}\label{spher}
We first consider spherical coordinates, 
\begin{equation}
	x_{1}=r\cos\phi\sin\theta,\quad x_{2}=r\sin\phi\sin\theta,\quad
	x_{3}=r\cos\theta,
	\label{spherical}
\end{equation}
where $r\in \mathbb{R}^{+}$, $\theta\in[0,\pi]$ and 
$\phi\in[0,2\pi[$.
The electric field has the components $E^{r}$, $E^{\theta}$ and 
$E^{\phi}$ in these coordinates. 

We concentrate here on the axisymmetric case where 
$E^{r}$, $E^{\theta}$ depend only on $r$, $\theta$ and $\omega$, and 
where $E^{\phi}(r,\theta,\omega)=0$. 
Let $G$ be the twist potential, with 
\cite{MS},
\begin{equation}
    G = \frac{1}{r}\left((rE^{\theta})_{r}-E^{r}_{\theta}\right)
    \label{NAsp1}.
\end{equation}
 For the Helmholtz equation (\ref{Helmholtz}) we get with \cite{MS} 
 and (\ref{NAsp1})
\begin{align}
(G\sin\theta)_{\theta}\frac{1}{\sin\theta}-\omega^{2}\epsilon rE^{r}  
&=rf^r, \nonumber\\
(rG)_{r}+\omega^{2}\epsilon r E^{\theta}  &=-rf^{\theta}
    \label{NAsp2}.
\end{align}
This system, together with some boundary conditions to be detailed in 
the following section, determines $E^{r}$ and $E^{\theta}$. 

Putting $Y=G\sin\theta$ and $x=\cos\theta$, system (\ref{NAsp2}) is 
equivalent to 
\begin{align}
Y_{x}+\omega^{2}\epsilon rE^{r}  
&=-rf^r, \nonumber\\
(rY)_{r}+\omega^{2}\epsilon r\sin\theta E^{\theta}  
&=-r\sin\theta f^{\theta}
    \label{NAsp2a}.
\end{align}
By differentiating (\ref{NAsp2a}), 
we get for constant $\epsilon$
\begin{equation}
    r^{2}(Y_{rr}+\omega^{2}\epsilon Y)+2rY_{r}+(1-x^{2})Y_{xx}=f,
    \label{NAsp3}
\end{equation}
where $f=r\sin\theta(-(rf^{\theta})_{r}+f^{r}_{\theta})$.
For given $Y$, the components $E^{r}$ and $E^{\theta}$ can be 
obtained from (\ref{NAsp2a}) for non-vanishing $\epsilon$ and 
$\omega$. Thus in the axisymmetric case, the solution of the 
vector Helmholtz equation (\ref{Helmholtz}) is equivalent to the solution of 
the scalar Helmholtz-type equation (\ref{NAsp3}). Note that this 
equation is singular for $r=0$ and $r\to\infty$, and for $x=\pm1$, 
the axis of symmetry. 

Equation (\ref{NAsp3}) can be separated (we only 
consider the homogeneous equation here) and solved in terms of 
spherical Bessel functions $Z_{l}(y)$, i.e., solutions to
$$y^{2}Z_{l}''(y)+2yZ_{l}'(y)+(y^{2}-l(l+1))Z_{l}(y)=0,$$
and derivatives of  
Legendre polynomials $P_{l}(x)$, 
$l\in \mathbb{N}$, see \cite{AS}. The general 
formal solution of the homogeneous equation can thus be written in the form 
\begin{equation}
	Y = (1-x^{2})\sum_{l=1}^{\infty}a_{l}P_{l}'(x)Z_{l}(kr)
	\label{Ysep},
\end{equation}
where $a_{l}$, $l\in \mathbb{N}$ are constants, and where $k = \omega 
\sqrt{\epsilon}$. The functions $Z_{l}(kr)$, also called cylinder 
functions are linear combinations of the spherical Bessel or Neumann 
functions, or of the Hankel functions, see \cite{AS}. Near the 
origin, only the Bessel functions are regular, near infinity the 
Sommerfeld condition determines the corresponding Hankel functions. 

If one is interested in the solution of equation (\ref{NAsp3}) for 
arbitrary $f$, it might be useful to introduce the scalar Green's function 
$\mathbf{g}$
such that
\begin{equation}
	Y = \mathbf{g}\otimes f,
	\label{Greensp}
\end{equation}
where $\otimes$ denotes the convolution in $x$ and $r$. Formally 
$\mathbf{g}$ can be obtained by solving (\ref{NAsp3}) with 
$f=\delta^{(2)}$ where $\delta^{(2)}$ is the two-dimensional 
delta-function. In a similar way one can define the Green's function 
of the system (\ref{NAsp2a}), 
\begin{equation}
	\begin{pmatrix}
		E^{r} \\
		E^{\theta}
	\end{pmatrix} = 
	\begin{pmatrix}
		\mathcal{G}^{rr} & \mathcal{G}^{rx} \\
		\mathcal{G}^{xr} & \mathcal{G}^{xx}
	\end{pmatrix} \otimes
	\begin{pmatrix}
		f^{r} \\
		f^{\theta}
	\end{pmatrix}
	\label{Greensp2}.
\end{equation}
The entries of the Green's 
function in (\ref{Greensp2}) follow for non-vanishing $\omega$ and 
$\epsilon$ for given $\mathcal{G}$ from 
(\ref{NAsp2a}),
\begin{align}
	\omega^{2} \epsilon \mathcal{G}^{rr}& = 
	-\delta^{(2)}+\frac{1}{r}\mathbf{g}_{x}\otimes \left[
	r(1-x^{2})\partial_{x}\right],
	\nonumber\\
	 \omega^{2} \epsilon \mathcal{G}^{rx}&= \frac{1}{r}\mathbf{g}_{x}\otimes \left[
	r\sqrt{1-x^{2}} \partial_{r}r\right],
	\nonumber\\
	\omega^{2} \epsilon \mathcal{G}^{xr}& 
	=\frac{1}{r\sqrt{1-x^{2}}}(r\mathbf{g})_{r}\otimes \left[
	r(1-x^{2})\partial_{x}\right],
	\nonumber\\
	\omega^{2} \epsilon \mathcal{G}^{xx}& = -\delta^{(2)}
	+\frac{1}{r\sqrt{1-x^{2}}}(r\mathbf{g})_{r}\otimes \left[
	r \sqrt{1-x^{2}}\partial_{r}r\right].
	\label{Greensp3}
\end{align}

\subsection{Prolate spheroidal coordinates}
\begin{figure}[htb!]
  \includegraphics[width=0.49\textwidth]{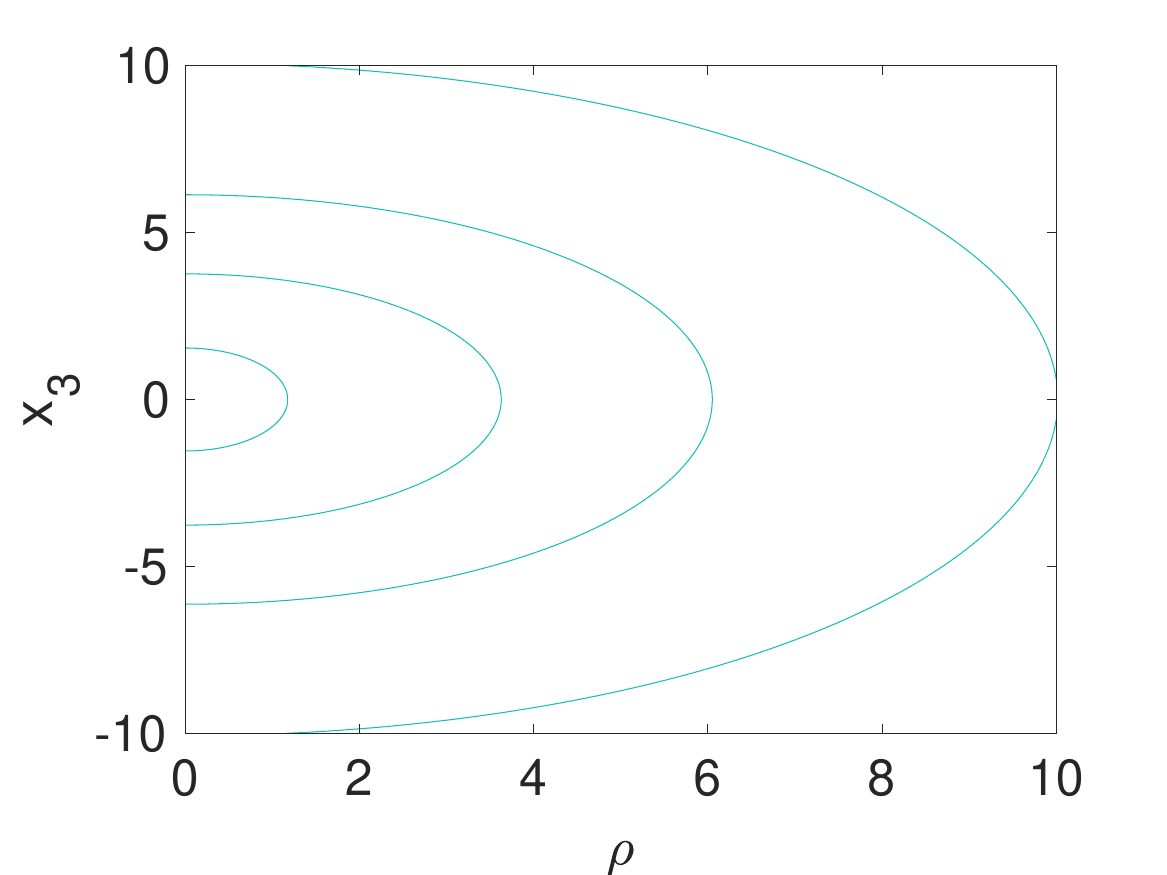}
  \includegraphics[width=0.49\textwidth]{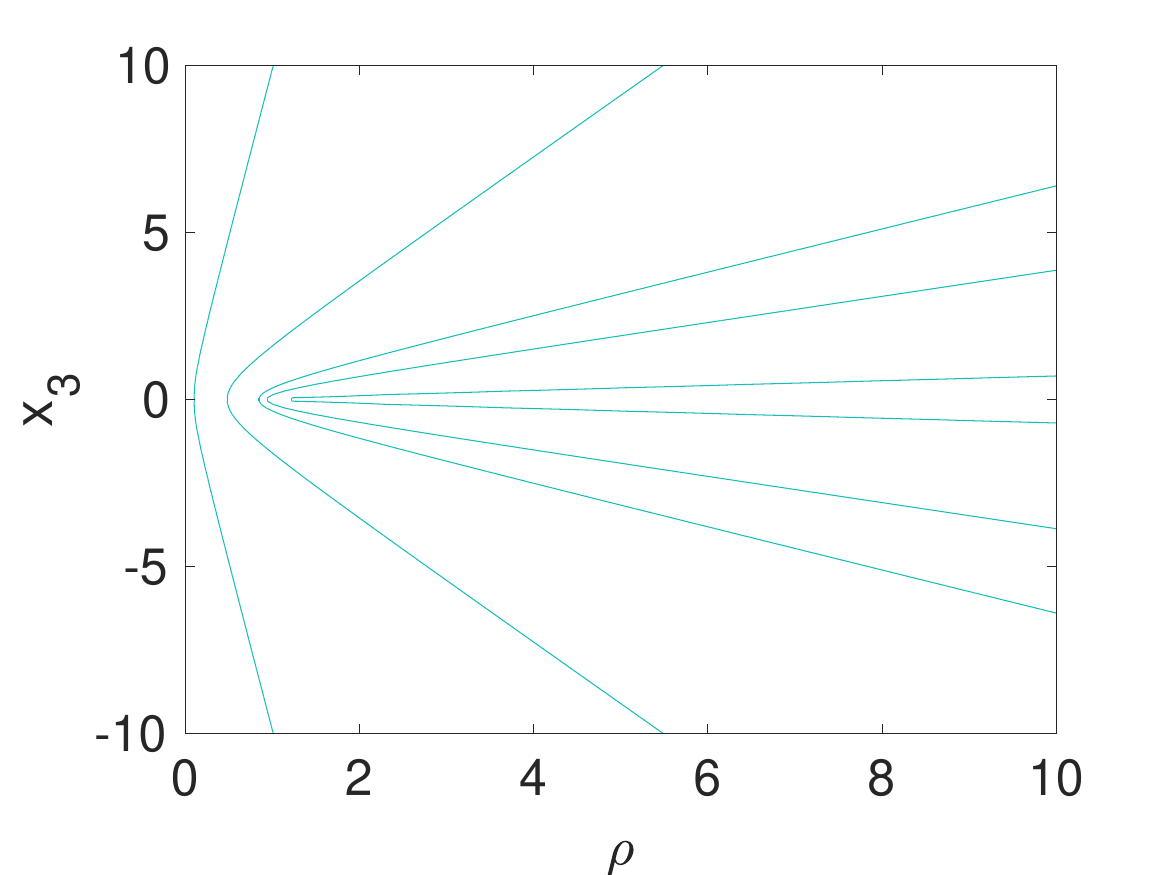}
\caption{Constant coordinate contours in the 
$\varrho,x_{3}$-plane, on the left (\ref{spheroids}) for $\eta=1,2,2.5,3$ 
(from left to right), and on the right (\ref{paraboloids}) for 
$\theta=0.1,0.5,1,1.2,1.5$ (also from left to right).}
 \label{CCS}
\end{figure}

Prolate spheroidal coordinates $\eta$, $\theta$ and $\phi$ 
with $0 \leq \eta < \infty$, $0 \leq \theta < \pi$ and $0 \leq \phi < 
2\pi$ are related to Cartesian coordinates via 
\begin{align}
x_{1}  &= a\sinh\eta\sin\theta\cos\phi , \nonumber\\
x_{2}  &= a\sinh\eta\sin\theta\sin\phi , \nonumber\\
x_{3}  &= a\cosh\eta\cos\theta \label{PS}. 
\end{align}
Constant coordinate surfaces are
\begin{equation}
    \frac{x_{1}^{2}+x_{2}^{2}}{a^{2}\sinh^{2}\eta}+\frac{x_{3}^{2}}{a^{2}\cosh^{2}\eta}=1,
    \label{spheroids}
\end{equation}
and 
\begin{equation}
    \frac{x_{1}^{2}+x_{2}^{2}}{a^{2}\sin^{2}\theta}-\frac{x_{3}^{2}}{a^{2}\cos^{2}\theta}=-1
    \label{paraboloids}.
\end{equation}
We show examples of these constant coordinate surfaces in 
Fig.~\ref{CCS} in the $\varrho,x_{3}$ plane where 
$\varrho:=\sqrt{x_{1}^{2}+x_{2}^{2}}$ in Fig.~\ref{CCS} (obviously there is a rotational 
symmetry with respect to the $\varrho=0$ axis). The spheroids can be 
seen on the left of the figure, the paraboloids on the right. 
The reader is referred to \cite{MS} for more information on these 
coordinates. 

As in the case of spherical coordinates, we concentrate on a 
situation with axial symmetry, i.e., $E^{\phi}=0$, and no dependence 
of $\mathbf{E}$
on the coordinate $\phi$. We introduce again the twist potential 
\begin{equation}
    F:=(\nabla\times 
    \mathbf{E})_{\phi}=\frac{1}{a\Psi}\left((\sqrt{\Psi}E^{\theta})_{\eta}-(\sqrt{\Psi}E^{\eta})_{\theta}\right)
    \label{NAps1},
\end{equation}
where we have put $\Psi:=\sinh^{2}\eta+\sin^{2}\theta$. 

This implies with \cite{MS} for the vector Helmholtz equation (\ref{Helmholtz})
\begin{align}
    (F\sin\theta)_{\theta}\frac{1}{\sin\theta}-\omega^{2}\epsilon a 
	\sqrt{\Psi}E^{\eta}
    & =a\sqrt{\Psi}f^{\eta},
    \nonumber\\
    (F\sinh\eta)_{\eta}\frac{1}{\sinh\eta}+\omega^{2}\epsilon a 
	\sqrt{\Psi}E^{\theta} &= {-}a\sqrt{\Psi} f^{\theta} 
    \label{NAps2a}.
\end{align}
The functions $E^{\eta}$ and $E^{\theta}$ can be determined from this 
system.

We put $x=\cos\theta$ and $y=\cosh\eta$,  which implies 
$\Psi=y^{2}-x^{2}$, and 
\begin{equation}
 X:=\sin\theta\sinh\eta F.
    \label{NAps2}
\end{equation}
The system (\ref{NAps2a}) then takes the form
\begin{align}
	X_{x}+\omega^{2}\epsilon a \sqrt{\psi}\sinh \eta E^{\eta} & 
	=-a\sqrt{\psi} \sinh\eta f^{\eta},
	\nonumber\\
	X_{y}+\omega^{2}\epsilon a \sqrt{\psi}\sin\theta E^{\theta} & =- 
	a\sqrt{\psi} \sin\theta f^{\theta}.
	\label{NAps2b}
\end{align}
By differentiating, the system (\ref{NAps2b}) is equivalent for 
constant $\epsilon$ to a scalar 
Helmholtz-type equation,
\begin{equation}
    (y^{2}-1)(X_{yy}+\omega^{2}\epsilon 
    a^{2}X)+(1-x^{2})(X_{xx}+\omega^{2}\epsilon 
    a^{2}X)=g
    \label{NAps3},
\end{equation}
where 
\begin{equation*}
g = -a\sqrt{(1-x^{2})(y^{2}-1)}\left[\sqrt{y^{2}-1}\left(f^\theta\sqrt{\Psi}
\right)_y 
+\sqrt{1-x^2}\left(f^\eta\sqrt{\Psi}\right)_x\right]. 
\end{equation*}

The homogeneous equation (\ref{NAps3}) can be separated in the form 
\begin{equation*}
X(y,x)=\sqrt{(1-x^{2})(y^{2}-1)}R_{\kappa}(y)S_{\kappa}(x),
\end{equation*}
where $\kappa$ is a complex constant, which leads to the ODEs 
\begin{align}
    (y^{2}-1)R_{\kappa}''(y)+2yR_{\kappa}'(y)+\left(\omega^{2}\epsilon 
    a^{2}(y^{2}-1)-\kappa-\frac{1}{y^{2}-1}\right)R_{\kappa}(y)& = 0,
    \nonumber\\
    (1-x^{2})S_{\kappa}''(x)-2xS_{\kappa}'(x)+\left(\omega^{2}\epsilon 
    a^{2}(1-x^{2})+\kappa-\frac{1}{1-x^{2}}\right)S_{\kappa}(x) & =0
    \label{NAps5}.
\end{align}
Its solutions are known as prolate spheroidal wave functions, see for 
instance \cite{AS}. The eigenvalues $\kappa$ are determined by the 
regularity of the solutions to the second equation on the axis 
($x=\pm1$). Thus the solution of (\ref{NAps3}) can be formally written 
in the form
\begin{equation}
	X = \sqrt{(1-x^{2})(y^{2}-1)}\sum_{\kappa\in K}^{}R_{\kappa}(y) 
	S_{\kappa}(x),
	\label{XODE}
\end{equation}
where $K$ denotes the (infinite) discrete spectrum of (\ref{NAps5}). 

If $\omega$ and $\epsilon$ do not vanish, $E^{\eta}$ and $E^{\theta}$ 
can be obtained via (\ref{NAps2}) for given $X$. Thus the scalar 
Helmholtz-type equation (\ref{NAps3}) is in this case equivalent to 
the Helmholtz equation (\ref{Helmholtz}). Equation (\ref{NAps3}) is 
singular for $y=1$ (the origin) and $y\to\infty$, and on the axis 
$x=\pm1$. 

As in the case of spherical coordinates in the previous subsection, 
it can  be useful to introduce the Green's function 
$\mathbf{g}$
such that
\begin{equation}
X = \mathbf{g}\otimes g
	\label{Greenps}
\end{equation}
where $\otimes$ denotes now the convolution in $x$ and $y$. The 
function  
$\mathbf{g}$ is the solution of  (\ref{NAps3}) with 
$g=\delta^{(2)}$ where $\delta^{(2)}$ is now the two-dimensional 
delta-function depending on $x$ and $y$. The Green's function 
of the system (\ref{NAps2b}) can be written in the form
\begin{equation}
	\begin{pmatrix}
		E^{\eta} \\
		E^{\theta}
	\end{pmatrix} = 
	\begin{pmatrix}
		\mathcal{G}^{yy} & \mathcal{G}^{yx} \\
		\mathcal{G}^{xy} & \mathcal{G}^{xx}
	\end{pmatrix}=
	\begin{pmatrix}
		f^{\eta} \\
		f^{\theta}
	\end{pmatrix}
	\label{Greenps2}.
\end{equation}

For non-vanishing $\omega$ and 
$\epsilon$ one gets for the Green's function in (\ref{Greenps2})
\begin{align}
	\omega^{2} \epsilon \mathcal{G}^{yy}& = 
	-\sqrt{y^{2}-1}\delta^{(2)}+\frac{1}{\sqrt{\psi}}\mathbf{g}_{x}\otimes \left[
	(y^{2}-1)\sqrt{1-x^{2}}\partial_{x}\sqrt{\psi}\right],
	\nonumber\\
	 \omega^{2} \epsilon \mathcal{G}^{yx}&= \frac{1}{\sqrt{\psi}}\mathbf{g}_{x}\otimes \left[
	(1-x^{2})\sqrt{y^{2}-1}\partial_{y}\sqrt{\psi}\right],
	\nonumber\\
	\omega^{2} \epsilon \mathcal{G}^{xy}& 
	=\frac{1}{\sqrt{\psi}}\mathbf{g}_{y}\otimes \left[
	(y^{2}-1)\sqrt{1-x^{2}}\partial_{x}\sqrt{\psi}\right],
	\nonumber\\
	\omega^{2} \epsilon \mathcal{G}^{xx}& = - \sqrt{1-x^{2}}\delta^{(2)}
	-\frac{1}{\sqrt{\psi}}\mathbf{g}_{y}\otimes \left[
	(1-x^{2})\sqrt{y^{2}-1}\partial_{y}\sqrt{\psi}\right].
	\label{Greenps3}
\end{align}

\section{Sommerfeld and matching conditions}\label{sec:bound}
In this section we summarize the matching conditions at the domain 
boundaries, and how the Sommerfeld radiation condition at infinity 
can be implemented. 

Since the Helmholtz equation is a second order PDE, 
one needs two matching conditions at each domain in order to obtain a 
unique solution. Since the domains are chosen in way that their 
boundary is a constant surface in  the `radial' coordinate, this 
means that a condition on the function and its normal derivative are 
required at each domain boundary. In domains where the equations are 
singular, the number of conditions can be less as detailed below.

\subsection{Sommerfeld radiation condition}
The Sommerfeld condition (\ref{som}) ensures that there is no 
incoming radiation from infinity. It implies that the solutions of 
the Helmholtz equation (\ref{Helmholtz}) can be written in the form 
\begin{equation}
	\mathbf{E} = e^{-i\omega ||\mathbf{x}||} 
	\mathbf{\tilde{E}}(||\mathbf{x}||,\omega)
	\label{ansatz},
\end{equation}
where $\mathbf{\tilde{E}}$ is a non-oscillatory function vanishing 
for $||\mathbf{x}||\to\infty$. Note that we assume that 
$\epsilon(||\mathbf{x}||,\omega)=1$ in an open environment of 
$\infty$, i.e., that all matter in the studied models is of finite volume. 

In spherical coordinates this implies that we can introduce in the 
vicinity of infinity the local parameter $\rho=1/r$ and split off the 
oscillatory terms as in (\ref{ansatz}), 
\begin{equation}
	E^{r} = e^{-i\omega r}\tilde{E}^{r},\quad 
	E^{\theta} = e^{-i\omega r}\tilde{E}^{\theta},\quad 
	Y = e^{-i\omega r}\tilde{Y}.
	\label{ansatz2}
\end{equation}

Thus we get for equation (\ref{NAsp3}) near infinity
\begin{equation}
	\rho^{2}\tilde{Y}_{\rho\rho}+2i\omega\tilde{Y}_{\rho}-\frac{2i\omega}{\rho}\tilde{Y}
	+(1-x^{2})\tilde{Y}_{xx}=e^{i\omega/\rho}f
	\label{Ytilde}.
\end{equation}
Note that this equation is singular for $\rho=0$ and $x=\pm1$, i.e., 
at infinity and on the symmetry axis. 

In the prolate spheroidal coordinates,  we make with (\ref{ansatz}) 
the ansatz
\begin{equation}
	E^{\eta} = e^{-i\omega ay}\tilde{E}^{\eta},\quad 
	E^{\theta} = e^{-i\omega ay}\tilde{E}^{\theta},\quad 
	X = e^{-i\omega ay}\tilde{X}
	\label{ansatz3}
\end{equation}
and introduce the local variable $\xi=1/y$ near infinity. Thus we 
get for (\ref{NAps3})
\begin{equation}
	(1-\xi^{2})(\xi^{2}\tilde{X}_{\xi\xi}+2(\xi-i\omega a)
	\tilde{X}_{\xi})+(1-x^{2})(\tilde{X}_{xx}+\omega^{2}a^{2}\tilde{X}) = e^{i\omega a/\xi}g
	\label{Xtilde}.
\end{equation}
Note that this equation is singular for $\xi=0$ and $x=\pm1$, i.e., 
at infinity and on the symmetry axis. 

\subsection{Matching conditions}\label{submatch}
The spectral methods we intend to apply in this paper are especially 
efficient if the physical boundaries coincide with domain boundaries, 
i.e., with constant coordinate surfaces. Therefore we discuss in this 
paper two sets of spheroidal coordinates and assume that $\epsilon$ 
is smooth or constant except for a finite number of values for the 
coordinate $r$ and $\eta$ respectively. 

It is known that in the absence of surface densities and currents, 
the normal components of $\mathbf{D}$ and $\mathbf{B}$ are continuous 
at the boundary, and that the same holds for the tangential 
components of $\mathbf{E}$ and $\mathbf{H}$. This means that in our 
cases
$E^{\theta}$ is continuous as well as $\epsilon E^{r}$ or $\epsilon 
E^{\eta}$ in spherical or prolate spheroidal coordinates 
respectively. 

The matching conditions for the function $Y$ can be read off from 
(\ref{NAsp2}): $Y$ and $Y_{r}/\epsilon $ are continuous at the 
boundaries. In a similar way the system (\ref{NAps2a}) gives the 
matching conditions for the function $X$: $X$ as well as $ 
X_{\eta}/\epsilon$ are continuous at the boundaries. 

We will work with  $N_{d}+1$ domains which are defined by the constant 
radii $r^{(i)}\leq r^{(i+1)}$, $i=1,\ldots,N_{d}$, in spherical coordinates and the constants 
$y^{(i)}\leq y^{(i+1)}$, $i=1,\ldots,N_{d}$, in prolate spheroidal coordinates:\\
$i=1$: $r<r^{(1)}$ ($y<y^{(1)}$): Near the origin, a singularity of the 
equations, special conditions need to be imposed to ensure a regular 
solution. In the spherical case,  the 
function $Y$ vanishes for $r=0$, since the spherical Bessel functions 
in (\ref{Ysep}) vanish there for $l>0$:
\begin{equation}
	Y^{(1)}(0,x,\omega) = 0.
	\label{Yr}
\end{equation}
In the prolate spheroidal case, $X$ must vanish at the origin in 
order to have a regular solution. 
In addition we impose that the functions $Y$, $X$ are continuous, 
\begin{equation}
	Y^{(1)}(r^{(1)},x,\omega)=Y^{(2)}(r^{(1)},x,\omega)
	\label{YI},
\end{equation}
where $Y^{(1)}$ is the function in domain I, and $Y^{(2)}$ is the 
function in domain II. Similarly we have
\begin{equation}
	X^{(1)}(y^{(1)},x,\omega)=X^{(2)}(y^{(1)},x,\omega)
	\label{XI}.
\end{equation}

$1<i<N_{d}-1$: $r^{(i-1)}<r<r^{(i)}$ ($y^{(i-1)}<y<y^{(i)}$): here we have to impose two 
conditions since there is no  
singularity in the radial coordinates. We impose continuity of 
the respective function at $r^{(i)}$, and a 
condition on the radial derivative at $r^{(i-1)}$:
\begin{align}
	Y^{(i)}(r^{(i)},x,\omega)& =Y^{(i+1)}(r^{(i)},x,\omega),
	\nonumber\\
	\frac{Y^{(i-1)}_{r}(r^{(i-1)},x,\omega)}{\epsilon^{(i-1)}(r^{(i-1)},\omega)}& 
	=\frac{Y^{(i)}_{r}(r^{(i-1)},x,\omega)}{\epsilon^{(i)}(r^{(i-1)},\omega)}
	\label{YII},
\end{align}
where $\epsilon^{(i)}$ are the 
values of $\epsilon$ in the respective domain. 

In prolate spheroidal coordinates we have the matching conditions
\begin{align}
	X^{(i)}(y^{(i)},x,\omega)& =X^{(i+1)}(y^{(i)},x,\omega),
	\nonumber\\
	\frac{X^{(i-1)}_{y}(y^{(i-1)},x,\omega)}{\epsilon^{(i-1)}(y^{(i-1)},\omega)}& 
	=\frac{X^{(i)}_{y}(y^{(i-1)},x,\omega)}{\epsilon^{(i)}(y^{(i-1)},\omega)}
	\label{YII}.
\end{align}

$i = N_{d}$:\\
In the domain bordering to the infinite one, the above conditions 
change as follows
\begin{align}
	Y^{(N_{d})}(r^{(N_{d})},x,\omega)& =e^{-i\omega 
	 r^{(N_{d})}}\tilde{Y}(1/r^{(N_{d})},x,\omega),
	\nonumber\\
	\frac{Y^{(N_{d}-1)}_{r}(r^{(N_{d}-1)},x,\omega)}{\epsilon^{(N_{d}-1)}(r^{(N_{d}-1)},\omega)}& 
	=\frac{Y^{N_{d}}_{r}(r^{(N_{d}-1)},x,\omega)}{\epsilon^{(N_{d}-1)}(r^{(N_{d}-1)},\omega)}
	\label{YIIa},
\end{align}
and 
\begin{align}
	 X^{(N_{d})}(y^{(N_{d})},x,\omega)& =e^{ -i\omega 
	 ay^{(N_{d})}}\tilde{X}(1/y^{(N_{d})},x,\omega),
	\nonumber\\
	\frac{X^{(N_{d}-1)}_{y}(y^{(N_{d}-1)},x,\omega)}{\epsilon^{(N_{d}-1)}(y^{(N_{d}-1)},\omega)}& 
	=\frac{X^{(N_{d})}_{y}(y^{(N_{d}-1)},x,\omega)}{\epsilon^{(N_{d})}(y^{(N_{d}-1)},\omega)}
	\label{XIIa},
\end{align}

$r>r^{(N_{d})}$ ($y>y^{(N_{d})}$): infinity is a singularity of the 
equations (\ref{Ytilde}) and (\ref{Xtilde}), but the vanishing of the 
respective solution at infinity has to be imposed in spherical 
coordinates (no condition is needed in the prolate spheroidal case). In addition we 
impose  
the matching condition on the radial derivative,
\begin{equation}
	\frac{Y^{(N_{d})}_{r}(r^{(N_{d})},x,\omega)}{\epsilon^{(N_{d})}(r^{(N_{d})},\omega)} 
	=e^{-i\omega 
	 r^{(N_{d})}}(\tilde{Y}_{\rho}(1/r^{(N_{d})},x,\omega) - i\omega 
	 \tilde{Y}(1/r^{(N_{d})},x,\omega)).
	\label{YIII}
\end{equation}

In prolate spheroidal coordinates we have
\begin{equation}
	\frac{X^{(N_{d})}_{y}(y^{(N_{d})},x,\omega)}{\epsilon^{(N_{d})}(y^{(N_{d})},\omega)} 
	=e^{ -i\omega a
	 y^{(N_{d})}}(\tilde{X}_{\xi}(1/y^{(N_{d})},x,\omega)  - i\omega a
	 \tilde{X}(1/^{(N_{d})},x,\omega)).
	\label{XIII}
\end{equation}

\section{Numerical approach}
In this section we briefly describe the numerical approach for the Helmholtz equations. In the angular coordinate, we 
always use a Chebyshev collocation method, in the `radial' 
coordinate, we consider several domains such that the line 
$\mathbb{R}^{+}$ is completely covered by these domains (infinity is 
simply a point on the grid). On each domain we use once more a 
Chebyshev collocation method. The matching conditions are imposed via 
a $\tau$-method \cite{tau}.

The essence of spectral methods is to approximate functions on a 
finite interval via functions being globally smooth on the considered 
interval. It is known that analytical functions are approximated by 
spectral methods with an error decreasing exponentially with the number of 
collocation points. Here we apply a Chebyshev collocation method, see 
\cite{trefethen} for details: the function to be approximated is 
sampled on the Chebyshev points $l_{n}=\cos(n\pi/N)$,  
$n=0,.1,\ldots,N$ with $N\in\mathbb{N}$. A function $u(l)$ is 
approximated on the interval $[-1,1]$ by the Lagrange polynomial 
$p(l)$ of degree $N$
passing through the collocation points, $p(l_{n})=u(l_{n})$, 
$n=0,\ldots,N$. The derivative of $u$ with respect to the argument is 
approximated via the derivative of the Lagrange polynomial which 
leads to the action of a \emph{Chebyshev differentiation matrix} $D$ 
on the vector $\mathrm{u}$ with components 
$u(l_{0}),\ldots,u(l_{N})$, i.e., $u'\approx D \mathrm{u}$. These 
Chebyshev differentiation matrices can be found for instance in 
\cite{trefethen,WR}. 

For the angular variable $x$, this method can be applied as described 
above since $x\in[-1,1]$. In the radial coordinate ($r$ or $y$ in the 
spherical and prolate spheroidal coordinates respectively), we 
introduce a number of domains which are chosen such that $\epsilon$ 
is smooth on each domain. Thus there will be a collection of radii 
$r^{(i)}<r^{(i+1)}$ (we only describe in the following this case since 
the treatment in $y$ is analogous), $i=1,\ldots,N_{d}$. Domain I 
is given by $r\leq r^{(1)}$ and thus contains the origin, domain 
$N_{d}+1$ is defined via $r>r^{(N_{d})}$ is infinite and will be 
compactified. 
Each of the intervals $[r^{(i)},r^{(i+1)}]$, $i=0,\ldots,N_{d}$ 
($r_{0}=0$) is mapped to the interval $[-1,1]$ via $r = 
(1+l)/2r^{(i+1)}+(1-l)r^{(i)}$, $l\in[-1,1]$. On the infinite interval, 
we apply the mapping  $r = 2/r_{N_{d}}/(1+l)$. On each interval the 
standard Chebyshev collocation points are introduced as well as the 
Chebyshev differentiation matrices. Since the domains II to $N_{d}$ are all identical from 
a mathematical point of view, we discuss in the following only the 
case of three domains. A generalization to a larger number $N_{d}>2$ is 
straight forward. 

This approach, i.e., discretization in both $r$ and $x$, allows to 
approximate the equations (\ref{NAsp3}) and 
(\ref{NAps3}) via a system of ordinary differential equations (ODEs) (this is a standard tensor grid). 
In the infinite domain, we 
discretize the equations (\ref{Ytilde}) respectively (\ref{Xtilde}). 
The matching conditions of the previous section are imposed via a 
$\tau$-method. This means that the equations in each domain 
corresponding to the radii $r_{i}$ and $r_{i+1}$ are replaced by the 
matching conditions in subsection \ref{submatch}. For a given right 
hand side $f$ (discretized in the same way), this leads to an 
equation of the Form $AY=f$ for some invertible matrix $A$ after 
discretization of the differentiation operators. This gives the 
wanted solution in each domain after solving the resulting linear 
system. This will be done in Matlab with the command 
`backslash', i.e., with essentially 
Gaussian elimination.
If one is interested in the Green's function, one simply has to 
replace the delta function $\delta^{(2)}$ by the identity in the 
considered vector space. 

Note that though the solution is only constructed on the collocation 
points
it can be obtained at all points in $\mathbb{R}^{2}$, with prescribed precision, via interpolation. 
An efficient and numerically stable way to do this is via barycentric 
interpolation, see \cite{barycentric} and references therein. 

A Chebyshev collocation method as presented above is equivalent to an 
expansion of a function $u(l)$ in terms of Chebyshev polynomials 
$T_{n}(l)$, $n\in\mathbb{N}$, where $T_{n}(l) = \cos(n\arccos(l))$. 
This means one approximates $u$ via
$$u(l)\approx \sum_{n=0}^{N}a_{n}T_{n}(l).$$ 
The Chebyshev coefficients $a_{n}$ are determined via
$$u(l_{n})=\sum_{m=0}^{N}a_{m}T_{m}(l_{n}),\quad n=0,\ldots,N.$$
This corresponds to a \emph{Fast Cosine Transform} which is related 
to the fast Fourier transform, see \cite{trefethen}, and is thus a 
very efficient way to compute the coefficients at a computational 
cost of $\mathcal{O}(N\ln N)$ operations.  

The decrease of the spectral coefficients  for a 
smooth function in both coordinates is expected to be exponential. This 
allows to allocate the numerical resolution in an efficient way. We 
first choose the domain boundaries according to the physical 
situation, i.e., discontinuities in the permittivity $\epsilon$ will 
be located at constant coordinate surfaces by assumption, and these 
will be chosen to be domain boundaries. Independently of this there 
will be always one domain containing the origin and one in the 
vicinity of infinity, both singularities of our equations. There 
will be always at least one domain in between these two domains, so the 
minimal number of domains will be three. \\
In each of these domains, we choose the resolution such that the 
spectral coefficients decrease to the aimed at accuracy, here 
essentially machine 
precision. Since the condition of the $N\times N$ Chebyshev differentiation 
matrix $D^{2}$ is of order $\mathcal{O}(N^{4})$, see 
\cite{trefethen}, one should aim in each domain at a small number of 
points which will not only lead to better conditioned matrices, but 
also to faster codes since the total differentiation matrix for all 
domains has a block structure, except for the matching conditions. \\
In 
general one can do an exploratory  low resolution run to estimate the 
optimal number of domains, the location of their boundaries and the 
resolution in each them. Since spectral methods are efficient for 
situations with simple geometry, this can be done by hand (this is 
always possible in our examples). If one were 
interested in an adaptive approach, one would as in \cite{CC} check 
the three (to avoid that coefficients vanish for symmetry reasons) 
spectral coefficients with the highest index for each domain 
and for each coordinate, and then vary the resolution until the largest 
of them in modulus is smaller than the aimed at accuracy.

\section{Examples}
In this section we study the performance of the presented codes for 
examples showing typical features of solutions to the Helmholtz 
equation (\ref{Helmholtz}). To construct examples in spherical and 
prolate spheroidal coordinates, we use what is jokingly called 
Synge's method in a general relativistic context: we make an ansatz 
for the solution and compute the right hand side of 
(\ref{Helmholtz}). This gives obviously an exact solution to the equation with 
this specific right hand side, 
which is then to be reproduced. Note that the goal of the first two 
examples is 
to provide interesting explicit test cases for the codes, not necessarily 
to study physically interesting situations for which no exact 
solutions are known. 
For simplicity we consider for these examples only vacuum, i.e., $\epsilon=1$ everywhere. 
An application to a typical problem in nano-optics 
is studied in the last subsection, the interaction of a radiating 
dipole with a spherical nano-particle, here silver, see e.g. 
\cite{DRJCCG} and references therein. In the nano-particle, the 
permittivity is considered to be given by a Drude model. No exact solution is 
known in this case with a piecewise constant permittivity, but we are 
able to numerically resolve this situation in a way that the spectral 
coefficients decrease to machine precision which can be seen as 
indicating the numerical accuracy. Note that the radii are chosen in 
all examples as in the physical problem for convenience. 

The codes we apply here are written in Matlab which 
is an interpreter language. Thus Matlab timings can depend strongly 
on the way of coding and how much precompiled code is actually used. 
Therefore the timings have to be taken with a grain of salt. But for 
practical applications it is of course useful to know at least the 
order of magnitude of time a run takes. The example of Fig.~\ref{figsolspher} 
takes roughly 2 seconds on a laptop, the example in 
Fig.~\ref{figsolspherom10} and in the last subsection roughly 80 
seconds. 

\subsection{Spherical coordinates}
In the case of spherical coordinates, we expect the function $Y$ 
to 
vanish for $r=0$ and for $r\gg 1$ to be 
oscillatory of the form (\ref{asym}) in order to satisfy the 
Sommerfeld condition.  

As an example for a function with this behavior we consider 
\begin{equation}
	Y = \frac{re^{-i\omega r}}{1+r^{2}+x^{2}}
	\label{ex1}.
\end{equation}
This implies with (\ref{Ytilde})
\begin{equation}
	\begin{split}
	f &= e^{-1i\omega r}\left(\frac{2r(1-2i\omega r)}{1+r^{2}+x^{2}}+
    \frac{2r^3(2i\omega r-5)}{(1+r^{2}+x^{2})^2} + 
	\frac{8r^5}{(1+r^{2}+x^{2})^3}\right.\\ 
	&-\left.
    \frac{2r(1-x^{2})}{(1+r^{2}+x^{2})^2} + 
	\frac{8rx^{2}(1-x^{2})}{(1r^{2}+x^{2})^3}\right);
\end{split}
	\label{fex}
\end{equation}
The source $f$ does not tend to zero at infinity, but this is not 
necessarily unphysical since we have $f=r^{2} \sin 
\theta(rf^{\theta}_{r}+f^{r}_{\theta})$. Thus the source, for 
instance a free charge density, is multiplied by a factor $r^{2}$. 

\begin{figure}[htb!]
  \includegraphics[width=0.32\textwidth]{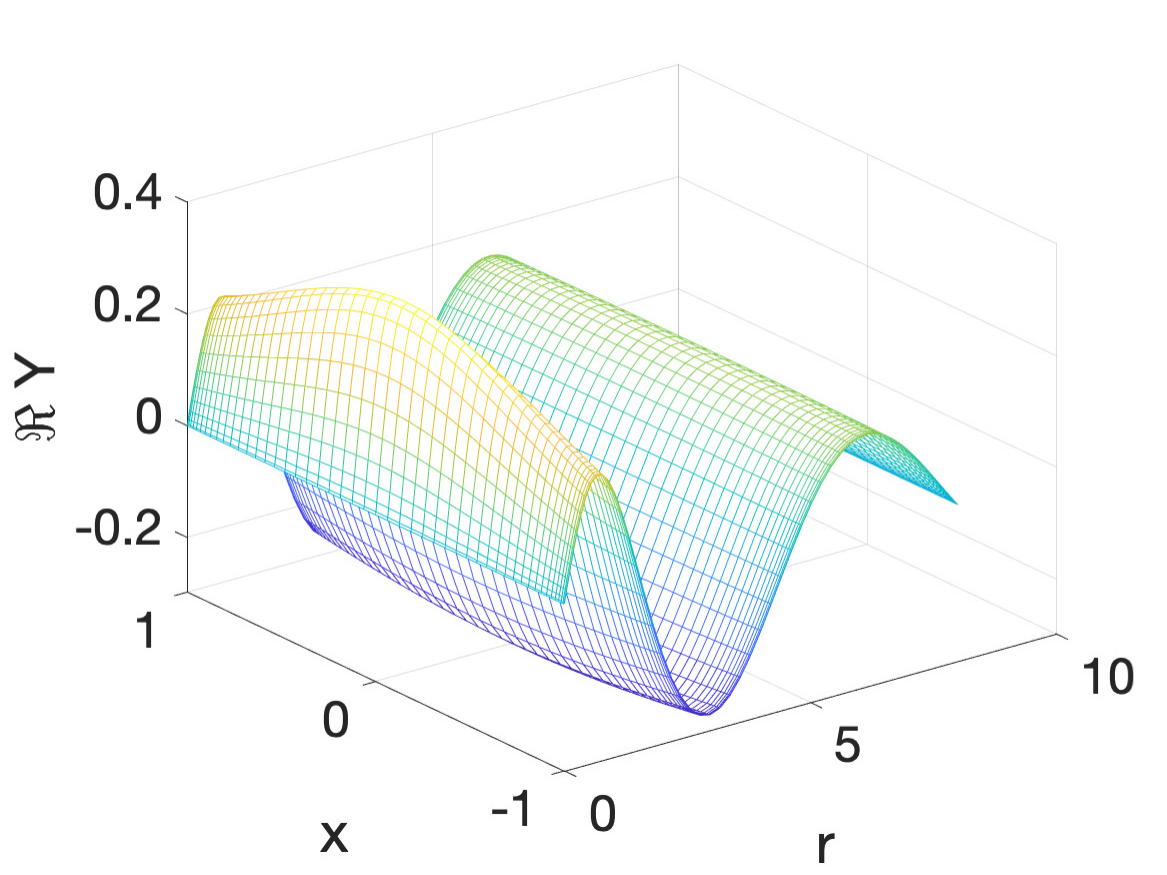}
  \includegraphics[width=0.32\textwidth]{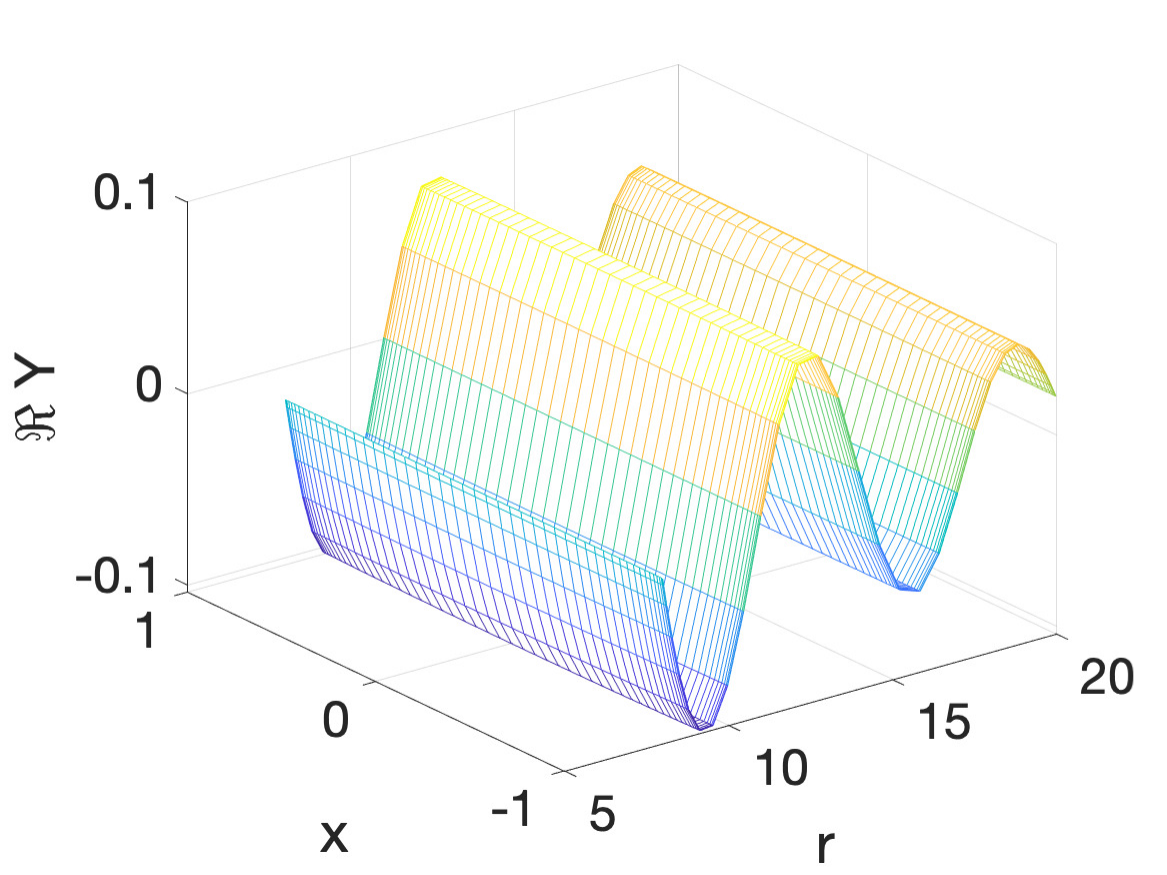}
  \includegraphics[width=0.32\textwidth]{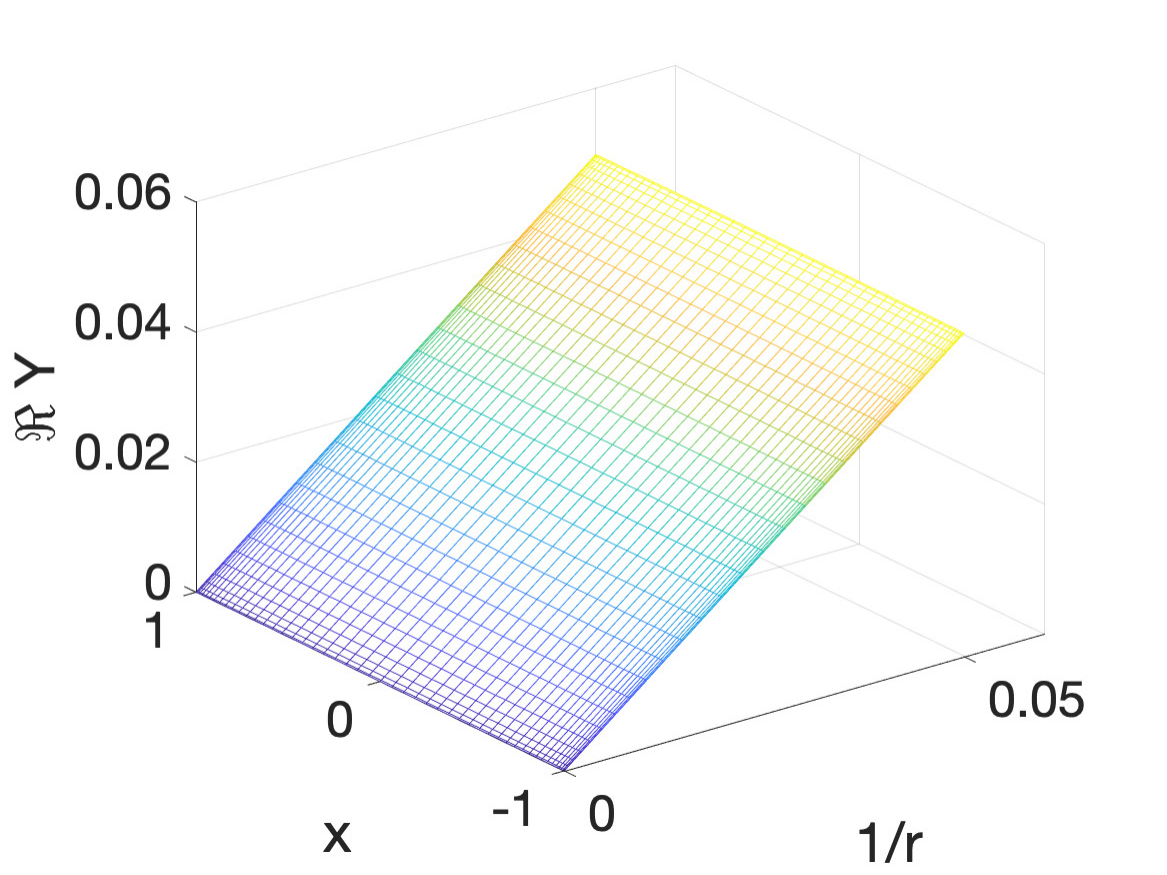}\\
   \includegraphics[width=0.32\textwidth]{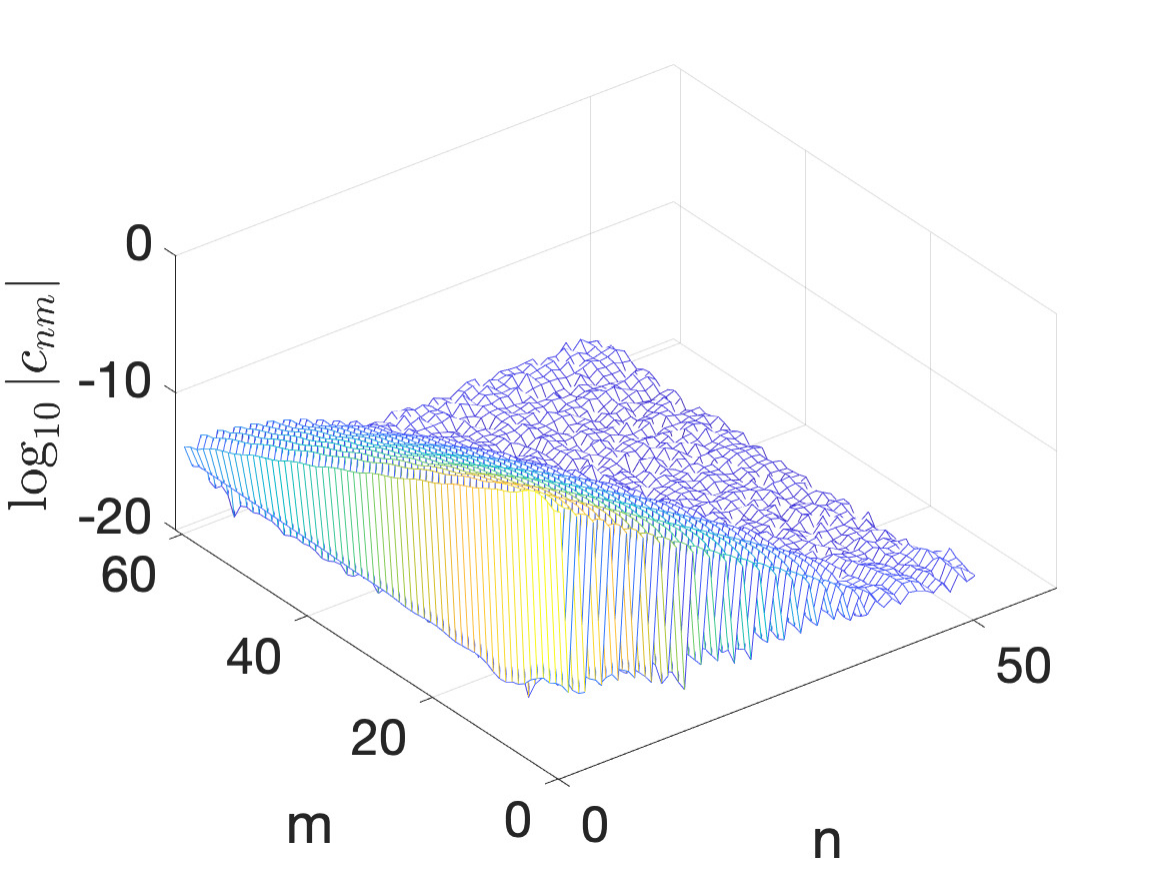}
  \includegraphics[width=0.32\textwidth]{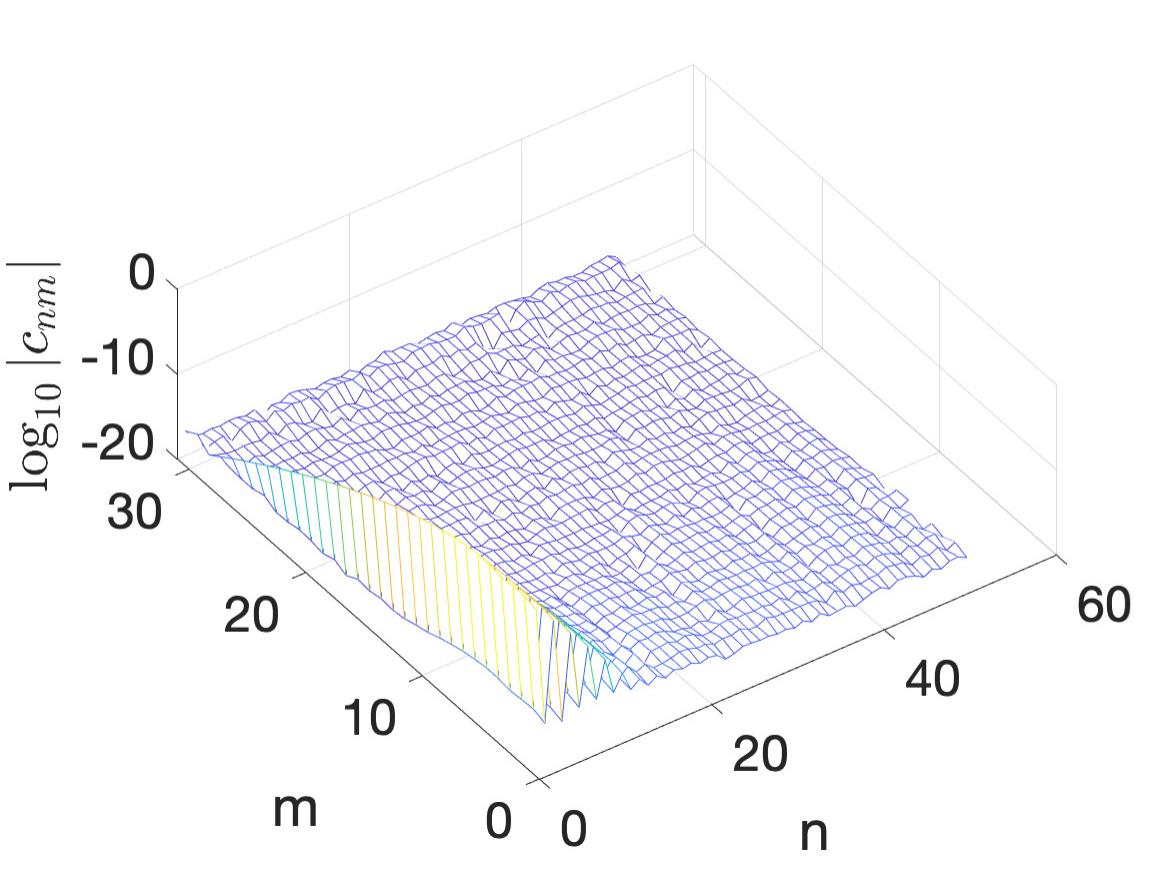}
  \includegraphics[width=0.32\textwidth]{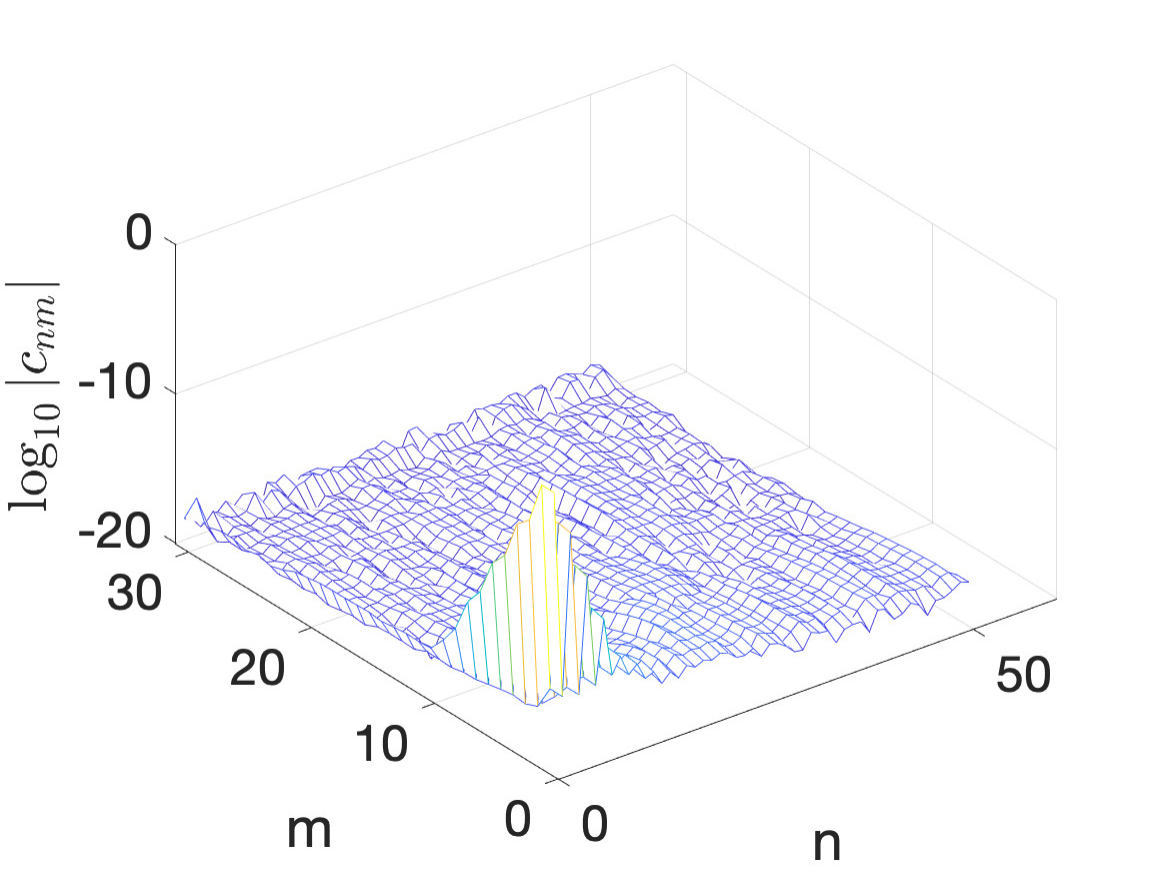}
\caption{Real part of the solution (\ref{ex1}) for $\omega=1$ in the 
 domains I, II, III (from left to right) in the upper row, and the 
 corresponding Chebyshev coefficients in the lower row.}
 \label{figsolspher}
\end{figure}

We use the three domains $r\leq 8$, $8<r<20$ and $r>20$. 
The real part of the solution is shown in these three domains for $\omega=1$ in 
the upper row of Fig.~\ref{figsolspher}. For the computation we use 
$N_{I}=60$, $N_{II}=30$, $N_{III}=30$ and $N_{x}=50$ Chebyshev 
polynomials. The Chebyshev coefficients in the respective domains can 
be seen in the lower row of Fig.~\ref{figsolspher}. It can be seen 
that they decrease with this choice of the number of collocation 
points to the order of machine precision. Note that the dependence of 
the solution (\ref{ex1}) on the variable $x$ is less pronounced the 
larger $r$ is. Thus one would be able to deal with less collocation 
points in domains II and III also in $x$, but in order to simplify 
the code, we use the same number of collocation points in $x$ in all 
domains. 

\begin{figure}[htb!]
  \includegraphics[width=0.32\textwidth]{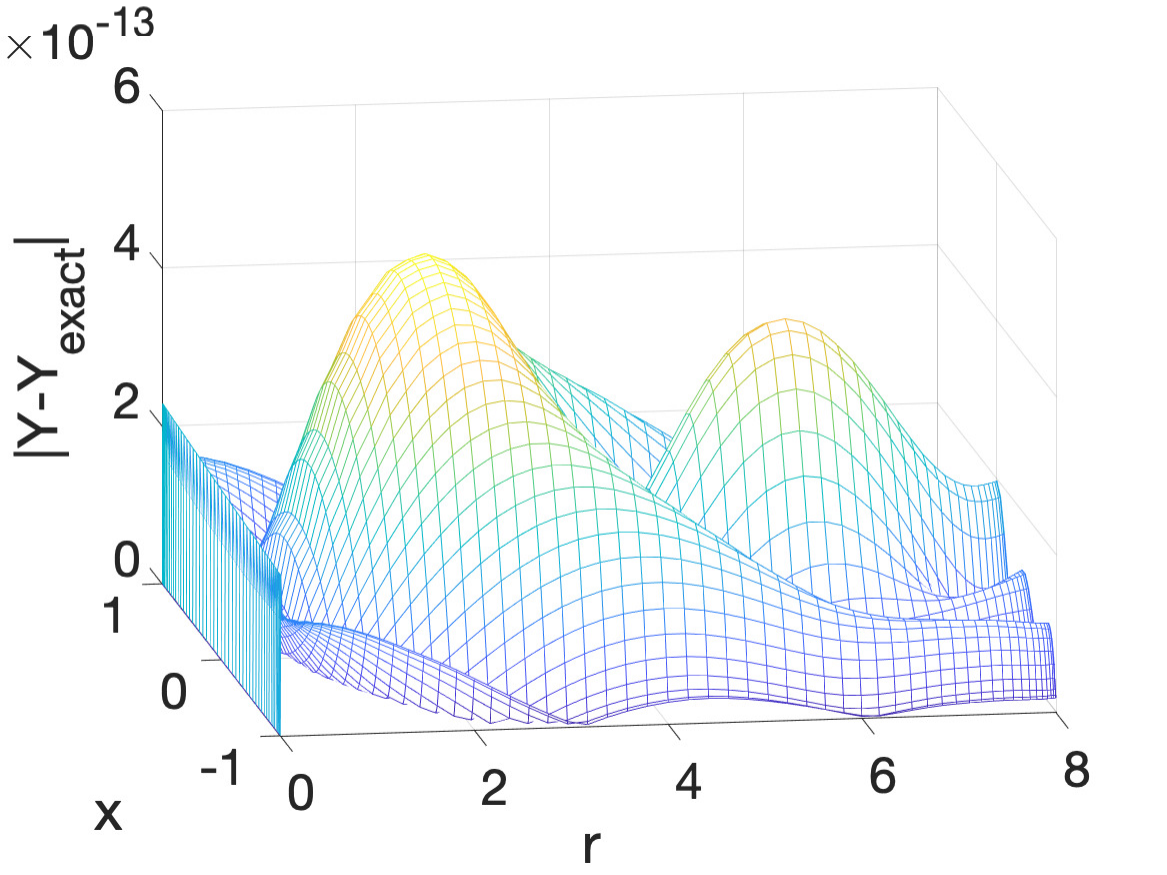}
  \includegraphics[width=0.32\textwidth]{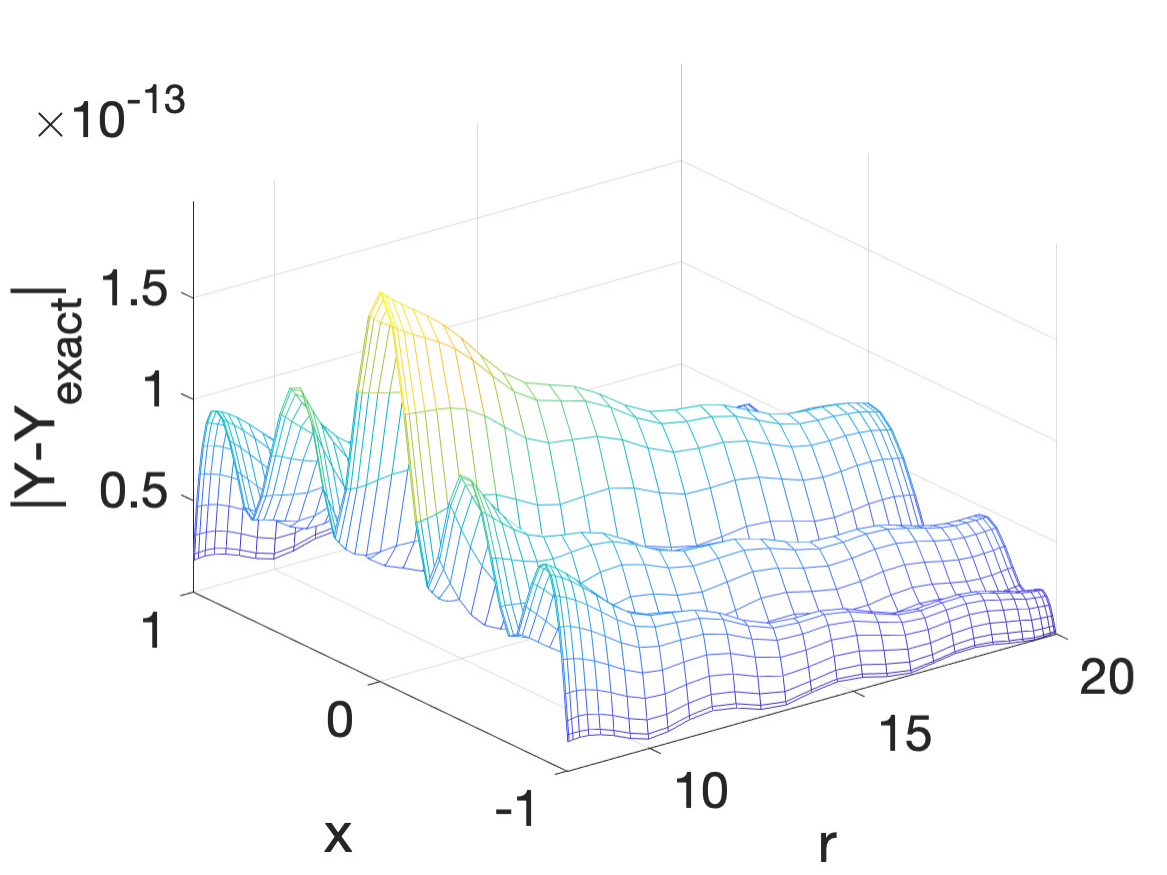}
  \includegraphics[width=0.32\textwidth]{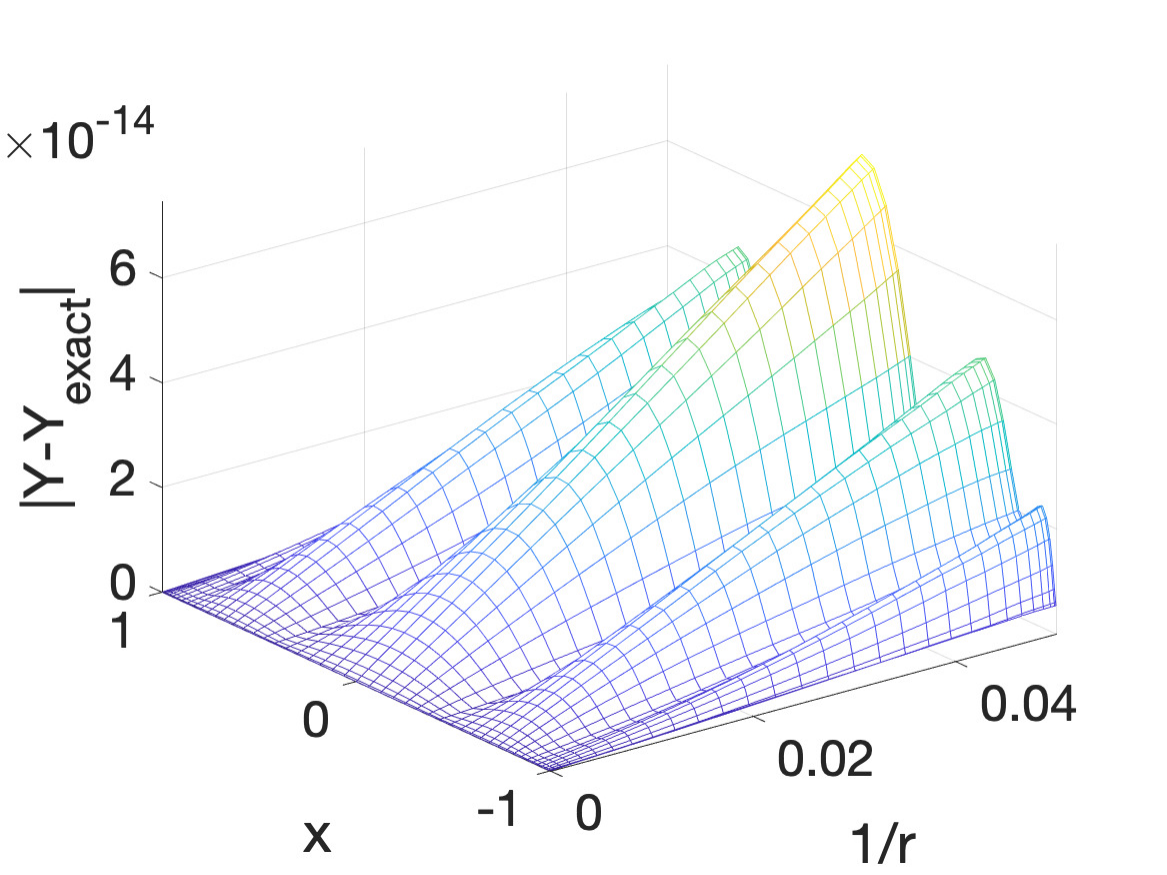}
\caption{Difference of the numerical solution of equation 
(\ref{Ytilde}) for the right hand side (\ref{fex}) and the exact 
solution (\ref{ex1}) in the 
 domains I, II, III (from left to right).}
 \label{figsolsphererr}
\end{figure}

If we solve equation (\ref{Ytilde}) for the right hand side 
(\ref{fex}) with the same number of collocation points as in 
Fig.~\ref{figsolspher}, one gets the difference between exact and 
numerical solution shown in Fig.~\ref{figsolsphererr}. It can be seen 
that it is globally of the order of $10^{-13}$, and thus as expected 
of the order as indicated by the highest Chebyshev coefficients in 
the lower row of Fig.~\ref{figsolspher}. 

\begin{figure}[htb!]
  \includegraphics[width=0.49\textwidth]{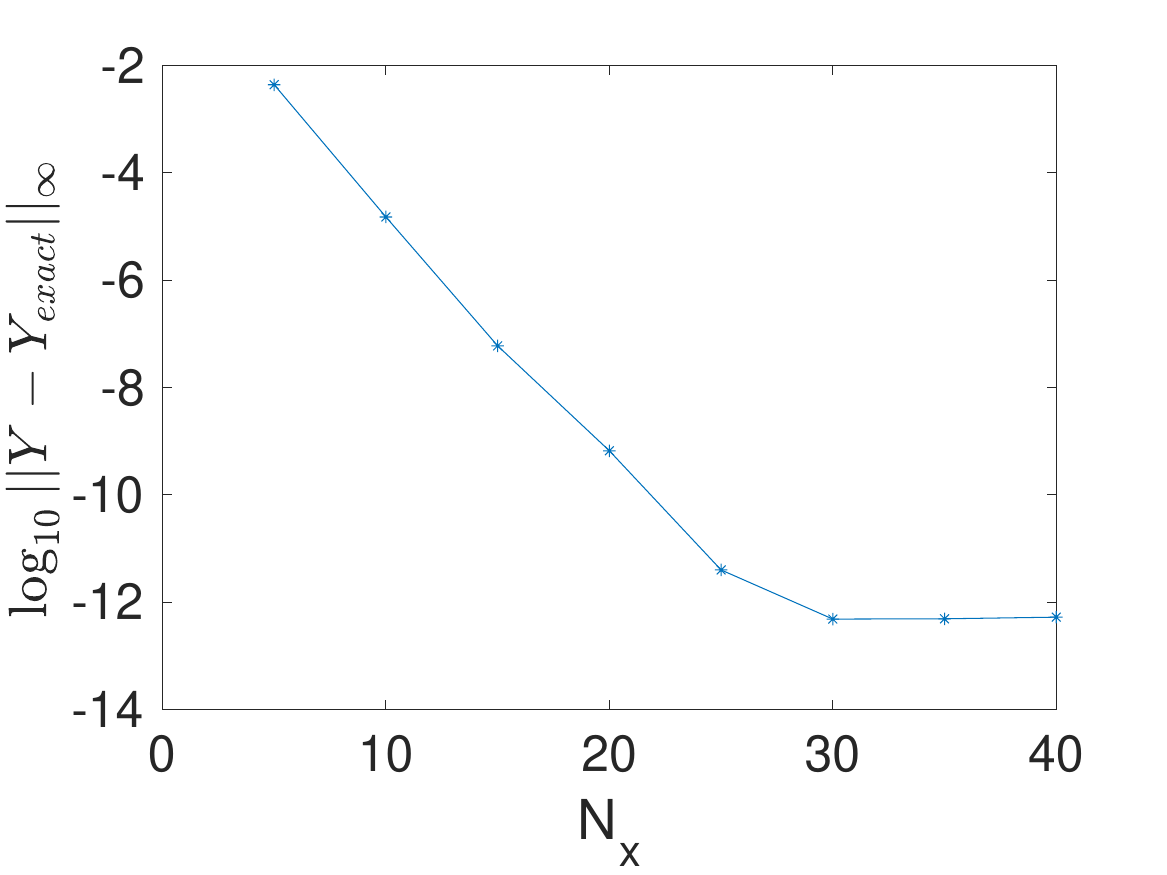}
  \includegraphics[width=0.49\textwidth]{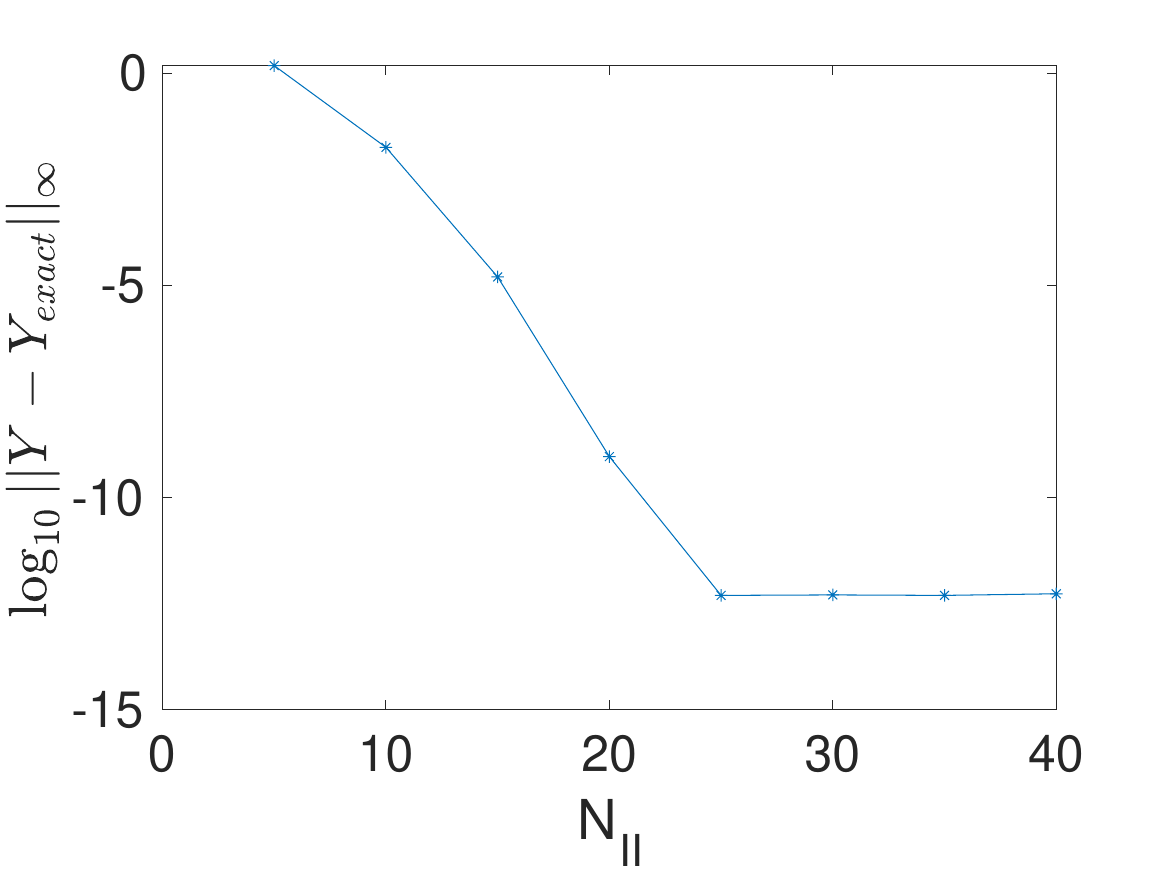}
\caption{$L^{\infty}$ norm of the difference of the numerical solution of equation 
(\ref{Ytilde}) for the right hand side (\ref{fex}) and the exact 
solution (\ref{ex1}) in dependence of $N_{x}$ on the left, and in 
dependence of $N_{II}$ on the right.}
 \label{figsolsphererrN}
\end{figure}

The dependence of the numerical error on the resolution in $x$ and 
$r$ can be seen in Fig.~\ref{figsolsphererrN}. For the same values of 
collocation points in $r$ as in Fig.~\ref{figsolspher}, the 
dependence of the difference between numerical and exact solution in 
the $L^{\infty}$ norm in dependence on $N_{x}$ can be seen on the 
left of Fig.~\ref{figsolsphererrN}. As expected it decreases 
exponentially and saturates essentially for $N_{x}\geq 30$. 

Note that though we use three domains in $r$, the solution in these 
domains are global for $r\in\mathbb{R}^{+}$. This is due to the fact 
that equation (\ref{Ytilde}) is elliptic, and that we impose on each 
domain boundary a $C^{1}$ condition on the solution. This leads to an 
analytical solution for $r\in\mathbb{R}^{+}$. Consequently a lack of 
resolution in one domain affects the numerical error in all domains. 
Thus to study the dependence of the numerical error on the resolution 
in $r$, it is sufficient to compute the global error in dependence of 
the resolution in just one domain. This error is shown for 
$N_{x}=50$, $N_{I}=60$ and $N_{III}=30$ in dependence of $N_{II}$ on 
the right of Fig.~\ref{figsolsphererrN}. The error decreases as 
expected exponentially with $N_{II}$ and saturates for $N_{II}\sim 
25$. 

\begin{figure}[htb!]
  \includegraphics[width=0.32\textwidth]{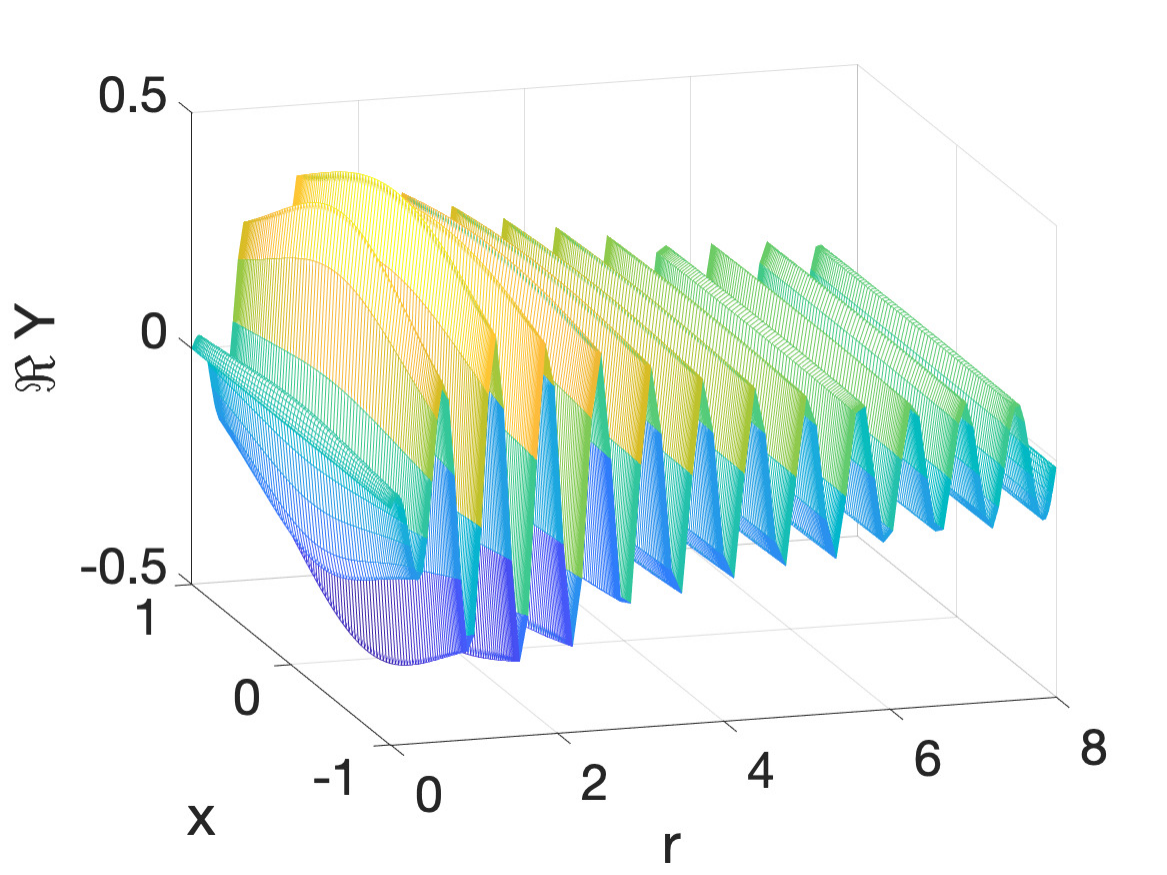}
  \includegraphics[width=0.32\textwidth]{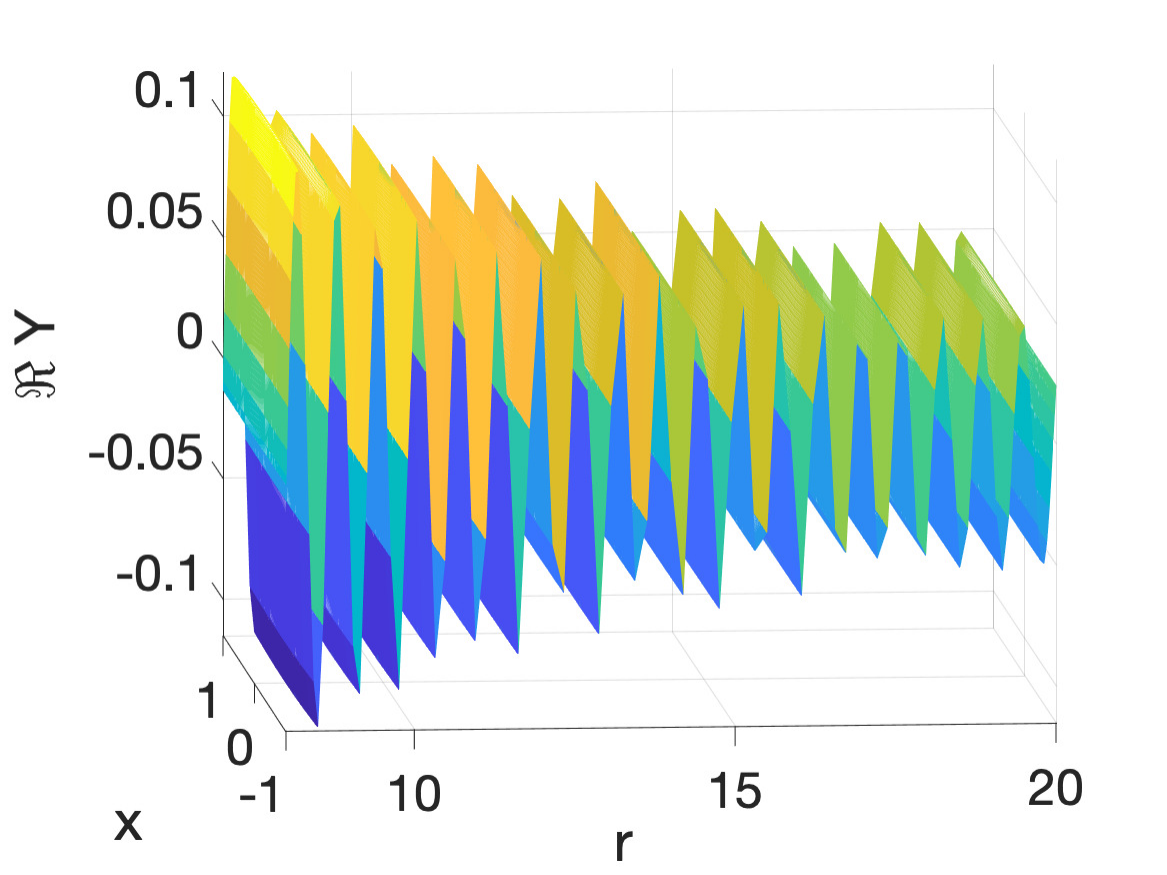}
  \includegraphics[width=0.32\textwidth]{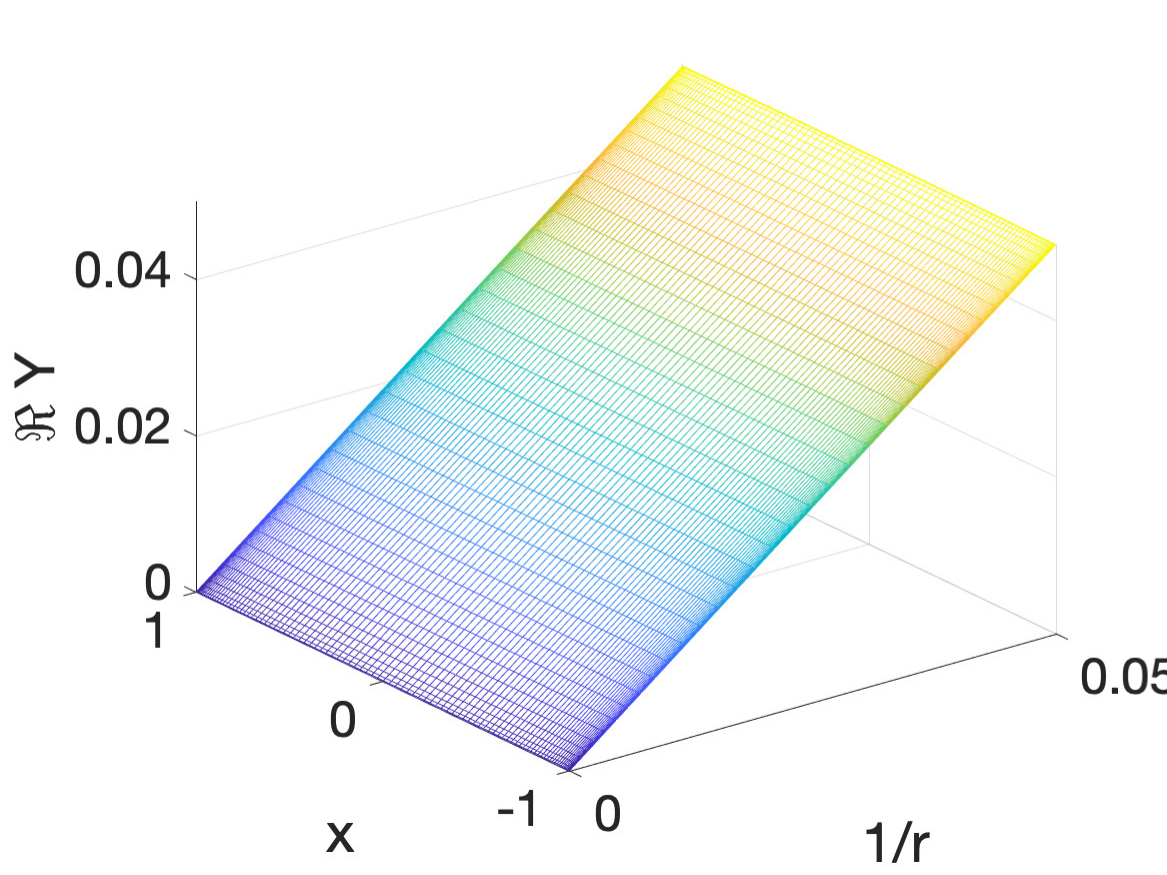}
\caption{Real part of the solution (\ref{ex1}) for $\omega=10$ in the 
 domains I, II, III (from left to right).}
 \label{figsolspherom10}
\end{figure}

Higher values of $\omega$ lead to a more oscillatory behavior of the 
solution, see Fig.~\ref{figsolspherom10} for $\omega=10$. This 
will make a higher resolution necessary. But with 
$N_{I}=N_{II}=100$, $N_{III}=40$ and $N_{x}=100$, 
we reproduce the solution (\ref{ex1})  to the order of 
$10^{-12}$ in this case. 

\subsection{Prolate spheroidal coordinates}
In the case of prolate spheroidal coordinates, we construct a similar 
test solution as in the spherical case: the solution must be
proportional to $\exp(-i\omega a 
y)/y$  for $y\to\infty$, and in addition must vanish for $y=1$. A 
possible candidate is
\begin{equation}
	X = e^{-i\omega a\sqrt{1+y^{2}}}\frac{y^{2}-1}{(y^{2}+x^{2})^{3/2}}
	\label{Xex}
\end{equation}
With (\ref{NAps3}) this implies that the source is of the form
\begin{equation}
	\begin{split}
	g &= (1-x^{2})(y^{2}-1)e^{-i\omega a\sqrt{1+y^{2}}}
	\left(\frac{3(4x^{2}-y^{2})}{(x^{2}+y^{2})^{7/2}}
	+\frac{\omega^{2}a^{2}}{(x^{2}+y^{2})^{3/2}}\right)\\
	&+(y^{2}-1)e^{-i\omega 
	a\sqrt{1+y^{2}}}\left(\frac{\omega^{2}a^{2}(y^{2}-1)}{(
	1+y^{2})(x^{2}+y^{2})^{3/2}}+\frac{2x^{2}+3-3y^{2}}{(x^{2}+y^{2})^{5/2}}
	-\frac{5y^{2}(2x^{2}+3-y^{2})}{(x^{2}+y^{2})^{7/2}}\right.\\
	&\left.-\frac{i\omega a y^{2}(2x^{2}+3-y^{2})}{\sqrt{1+y^{2}}
	(x^{2}+y^{2})^{5/2}}-\frac{i\omega 
	a(3y^{2}-1)}{\sqrt{1+y^{2}}(x^{2}+y^{2})^{3/2}}
	+\frac{i\omega 
	ay^{2}(y^{2}-1)(3+x^{2}+4y^{2})}{(1+y^{2})^{3/2}(x^{2}+y^{2})^{5/2}}
	\right).
	\end{split}
	\label{gex}
\end{equation}

\begin{figure}[htb!]
  \includegraphics[width=0.32\textwidth]{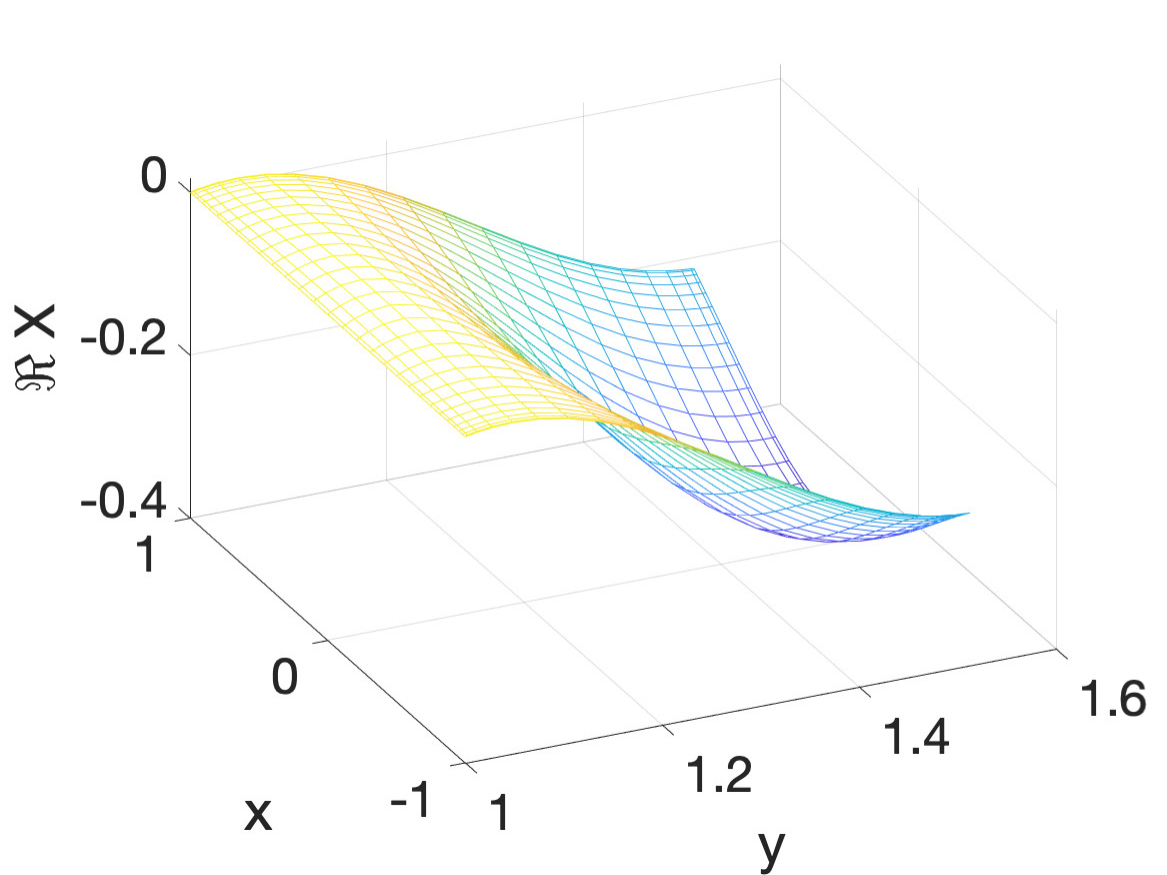}
  \includegraphics[width=0.32\textwidth]{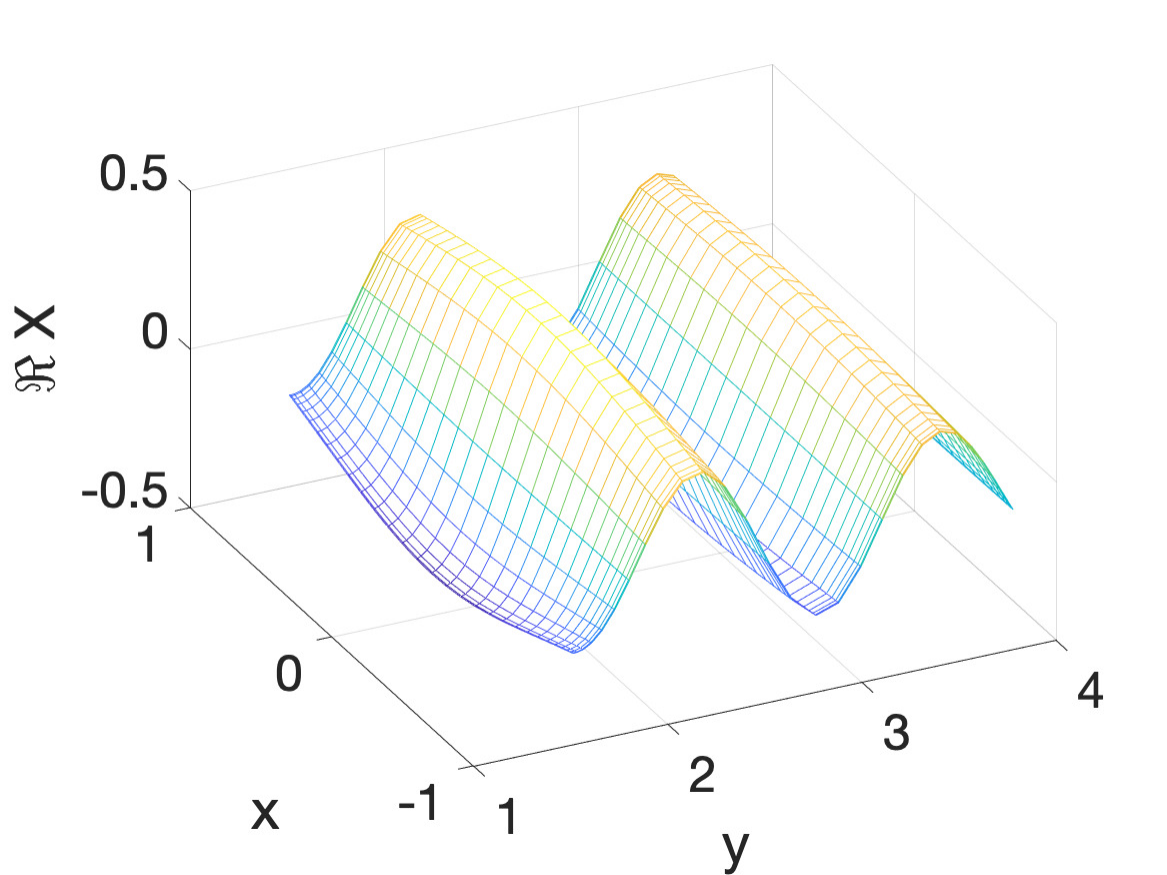}
  \includegraphics[width=0.32\textwidth]{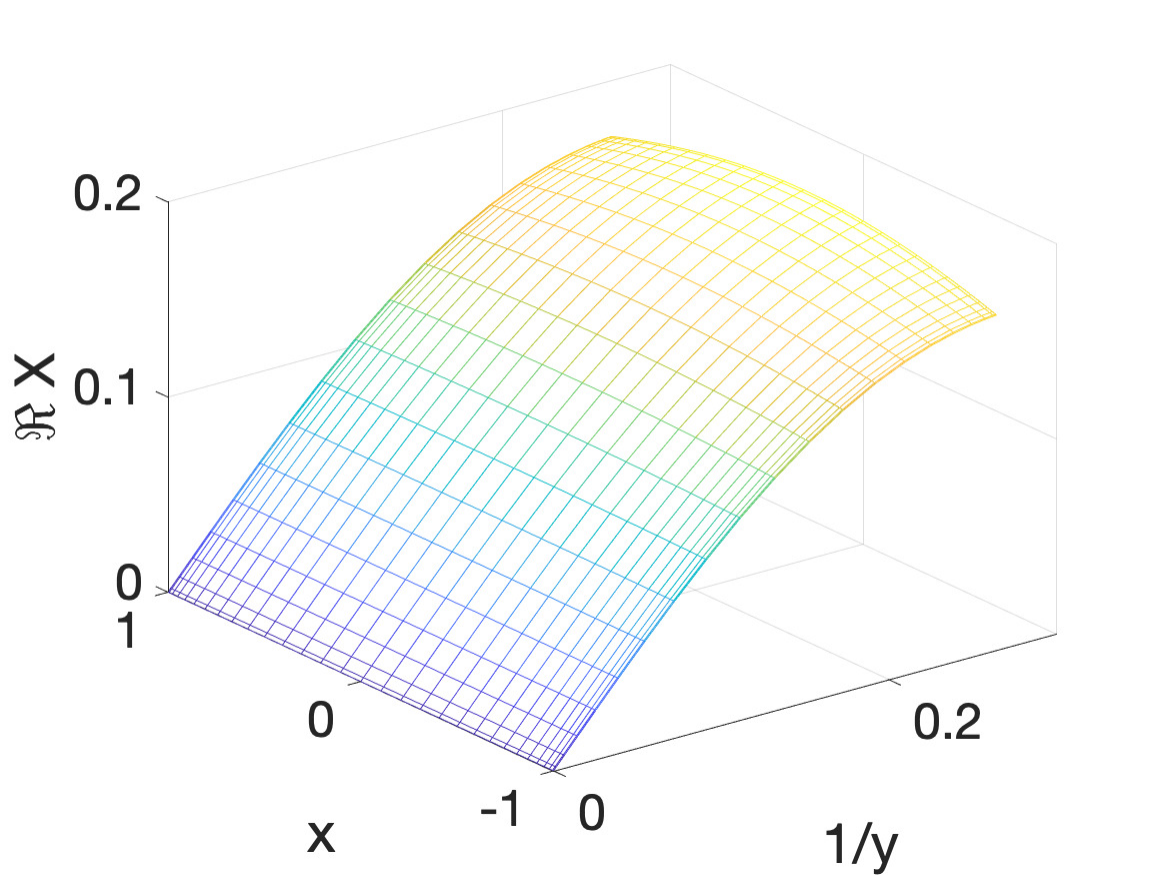}\\
   \includegraphics[width=0.32\textwidth]{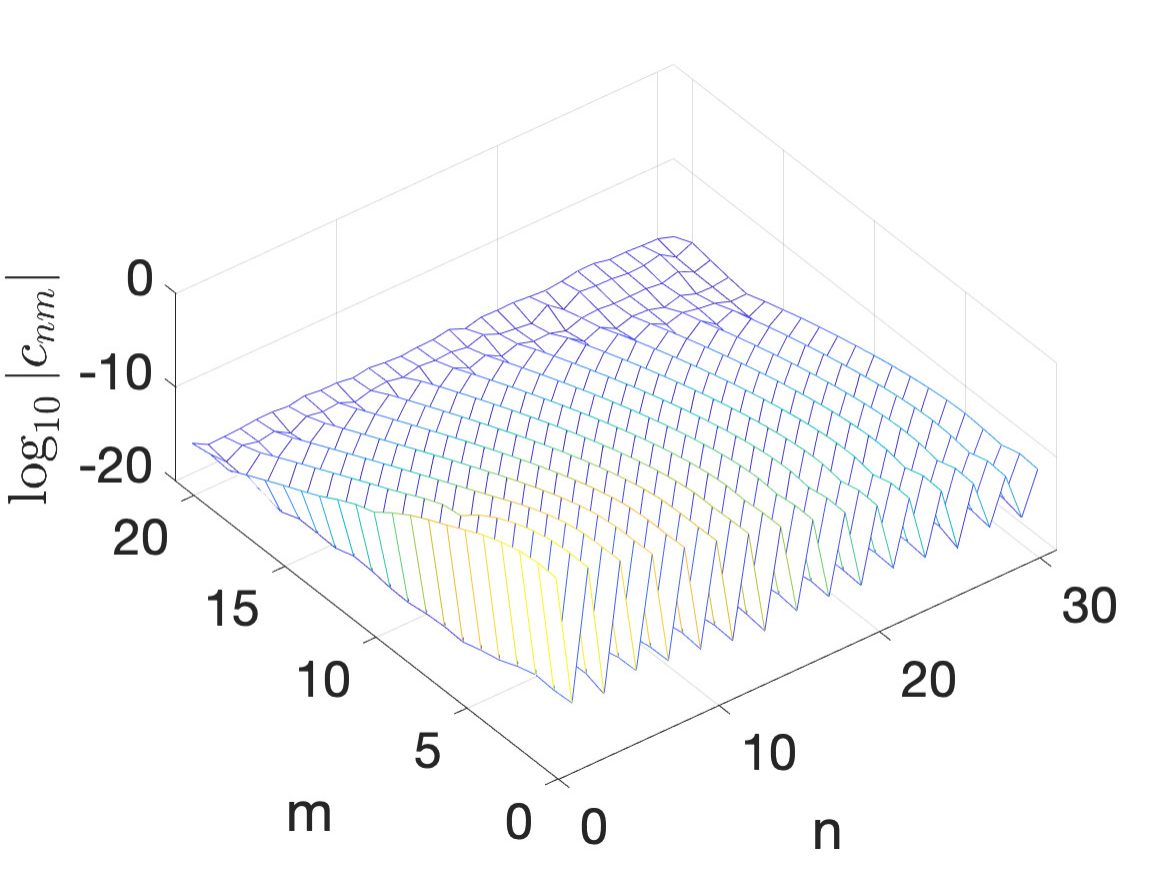}
  \includegraphics[width=0.32\textwidth]{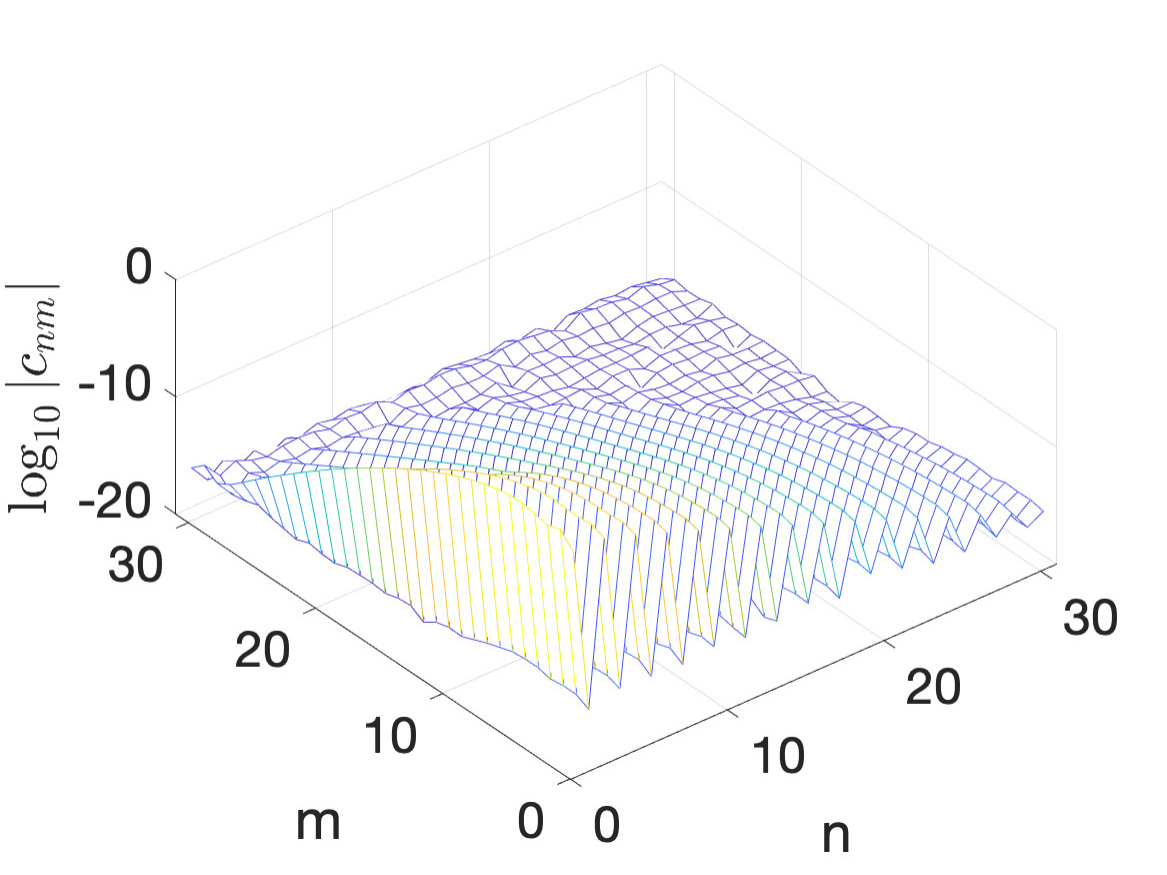}
  \includegraphics[width=0.32\textwidth]{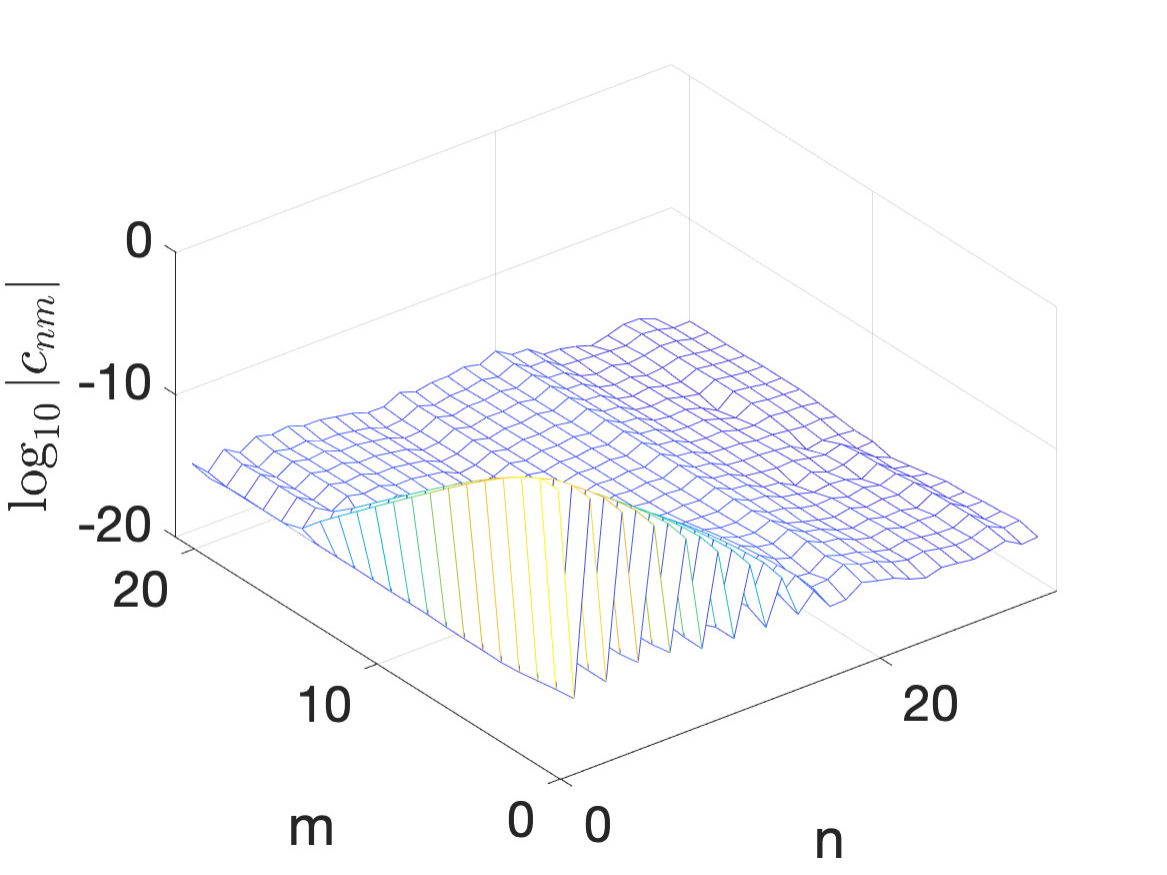}
\caption{Real part of the solution (\ref{Xex}) for $\omega=1$ in the 
 domains I, II, III (from left to right) in the upper row, and the 
 corresponding Chebyshev coefficients in the lower row.}
 \label{figsolps}
\end{figure}

We choose the domains with the same values on the axis as before, $a 
y^{I}=8$ and $ay^{II}=20$. If the boundary of domain I in the 
$x_{3}=0$ plane is $x_{1}=6$, one has $a \sim 5.29$. Solution 
(\ref{Xex}) for these values and $\omega=1$ can be seen in 
Fig. \ref{figsolps}. If we use $N_{I}=N_{III}=20$ and  
$N_{II}=N_{x}=30$ collocation points, we get the Chebyshev 
coefficients shown in the lower row of Fig.~\ref{figsolps}. They 
decrease in all cases to machine precision. 

\begin{figure}[htb!]
  \includegraphics[width=0.32\textwidth]{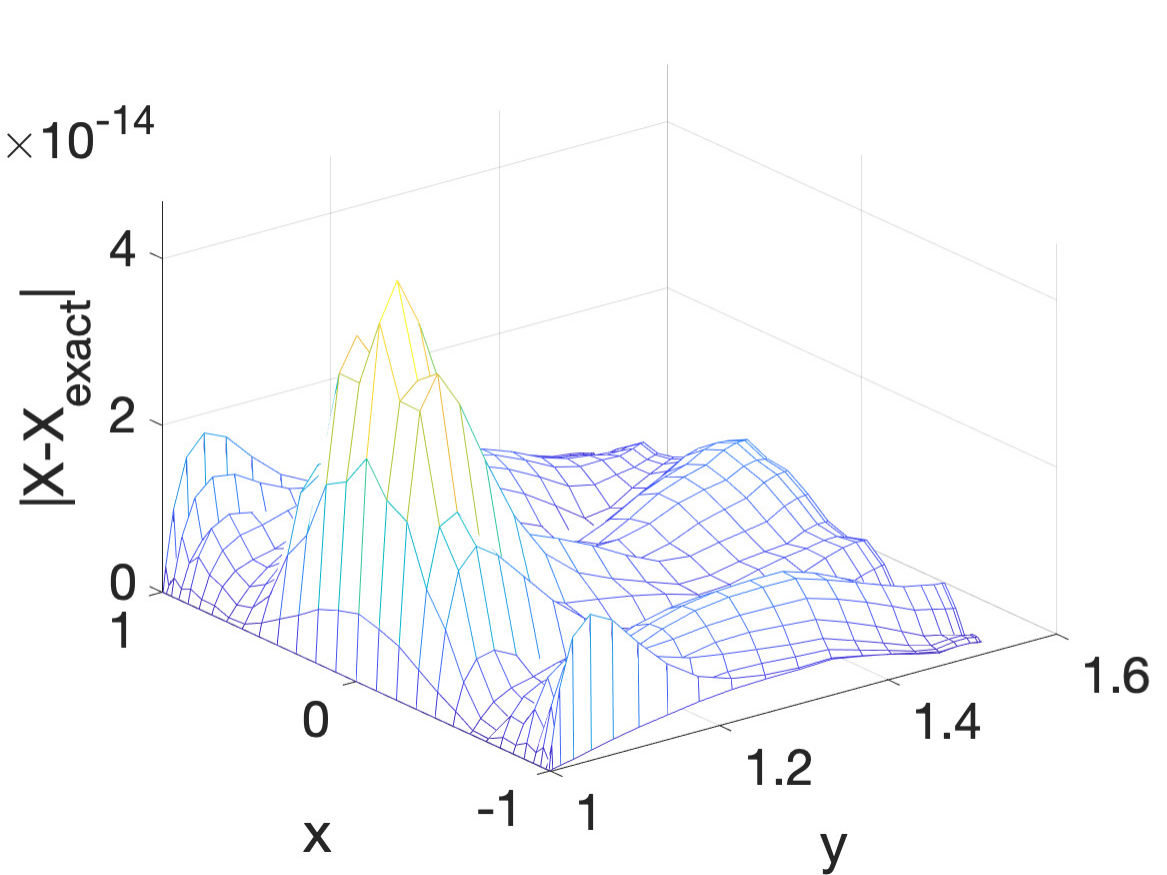}
  \includegraphics[width=0.32\textwidth]{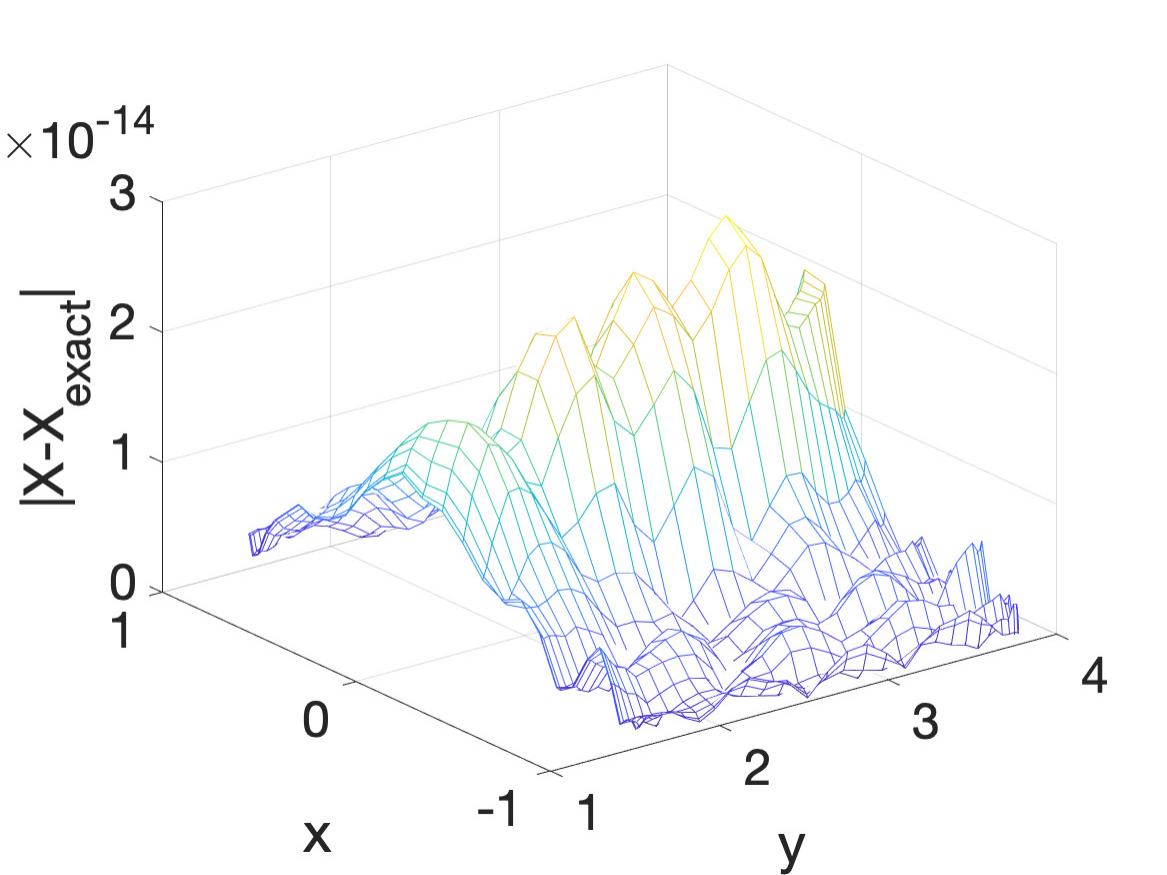}
  \includegraphics[width=0.32\textwidth]{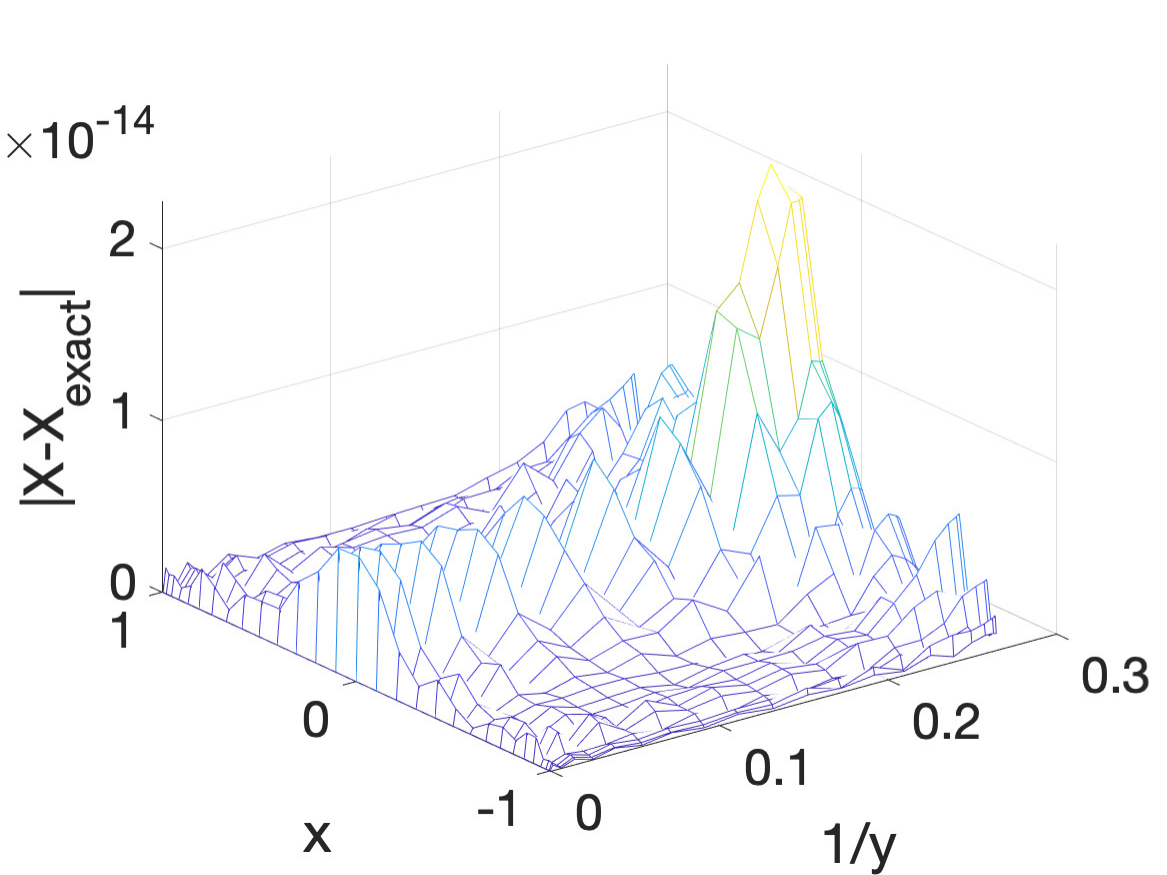}
\caption{Difference of the numerical solution of equation 
(\ref{Xtilde}) for the right hand side (\ref{gex}) and the exact 
solution (\ref{Xex}) in the 
 domains I, II, III (from left to right).}
 \label{figsolpserr}
\end{figure}

For this choice of the numerical parameters, the difference between 
numerical and exact solution is shown in Fig.~\ref{figsolpserr}. As 
expected from the Chebyshev coefficients, the error is globally of 
the order of $10^{-14}$ in this case. 

\begin{figure}[htb!]
  \includegraphics[width=0.49\textwidth]{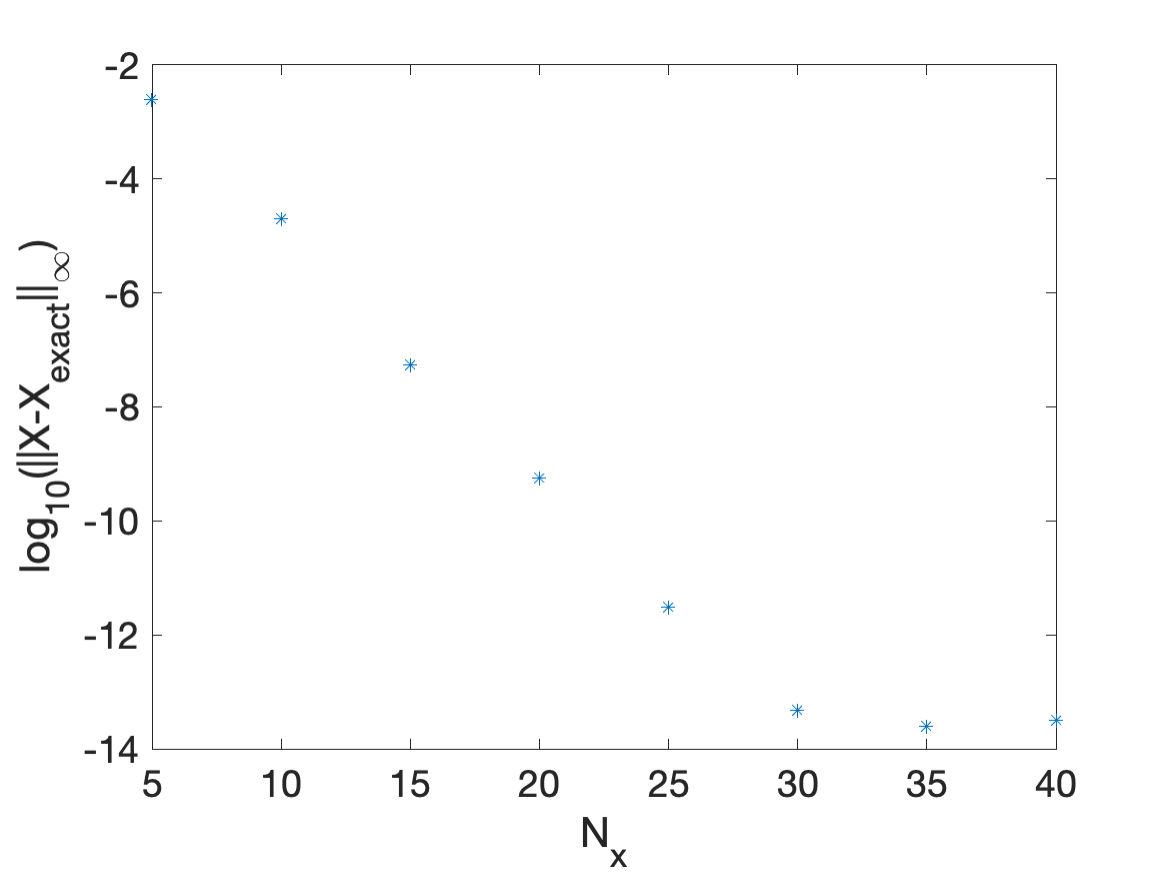}
  \includegraphics[width=0.49\textwidth]{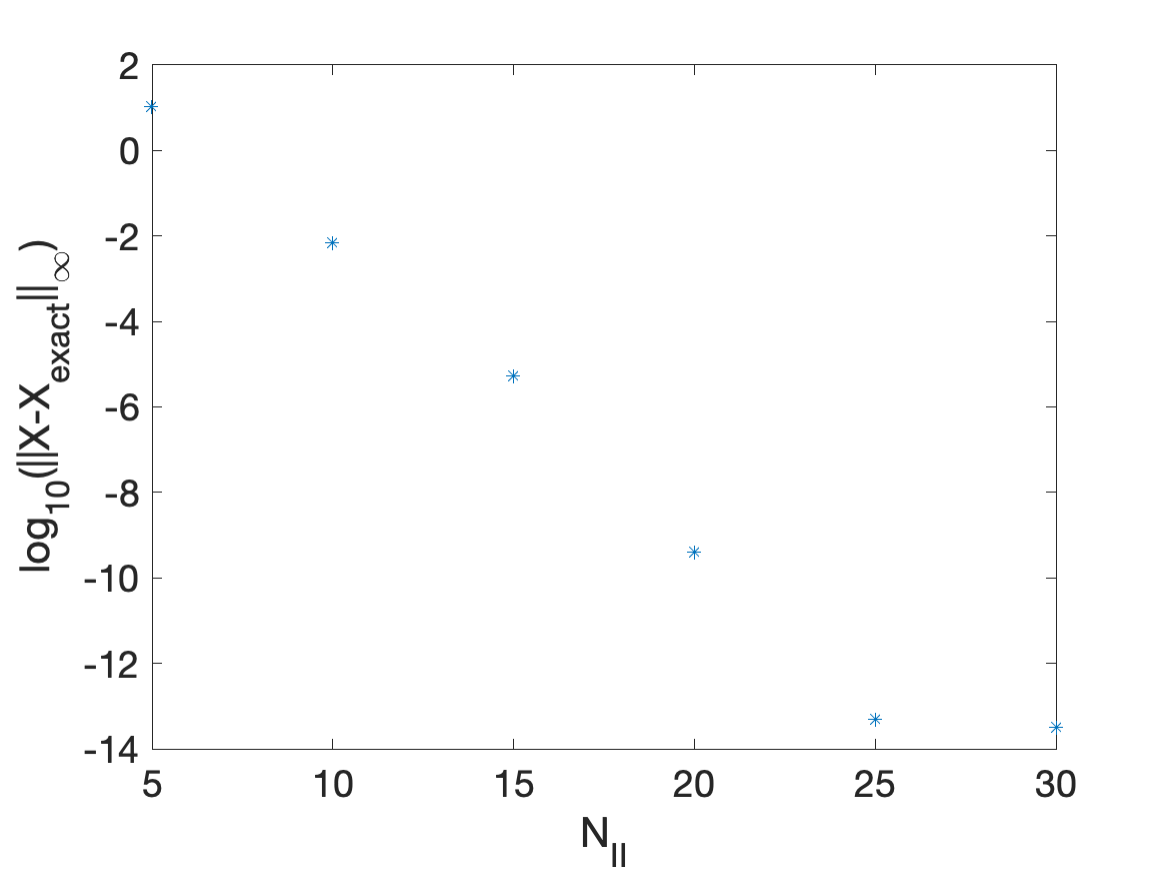}
\caption{$L^{\infty}$ norm of the 
difference of the numerical solution of equation 
(\ref{Xtilde}) for the right hand side (\ref{gex}) and the exact 
solution (\ref{Xex}) in dependence of $N_{x}$ on the left, and in 
dependence of $N_{II}$ on the right).}
 \label{figsolpserrN}
\end{figure}

The dependence of the numerical error on the resolution can be again 
studied by varying the number of collocation points. On the left of 
Fig.~\ref{figsolpserrN}, the same number of collocation points in $y$ 
is applied, and $N_{x}$ varies. It can be seen that the error 
decreases exponentially and saturates for $N_{x}\sim 30$. On the 
right of the same figure, $N_{II}$ is varied. Again the error 
decreases exponentially and saturates for $N_{II}\sim 25$. 

\begin{figure}[htb!]
  \includegraphics[width=0.56\textwidth]{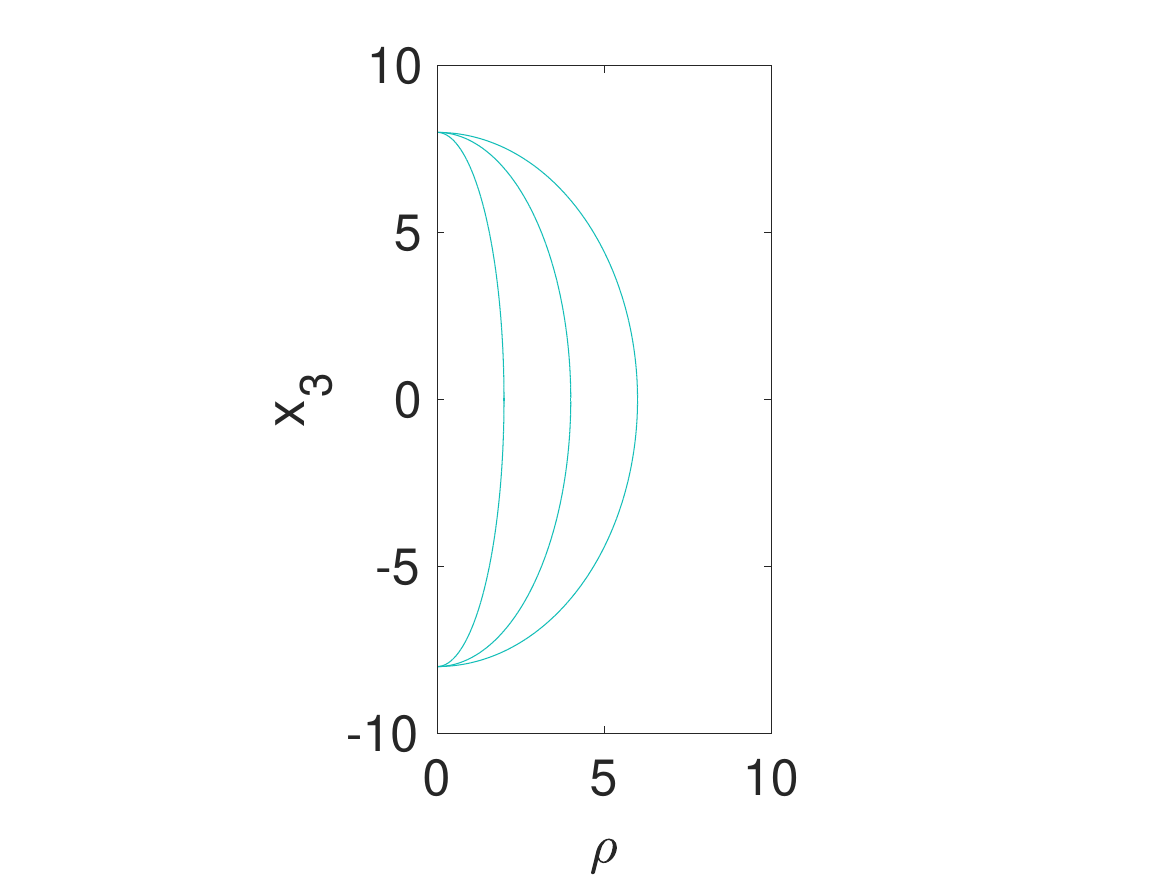}
  \includegraphics[width=0.43\textwidth]{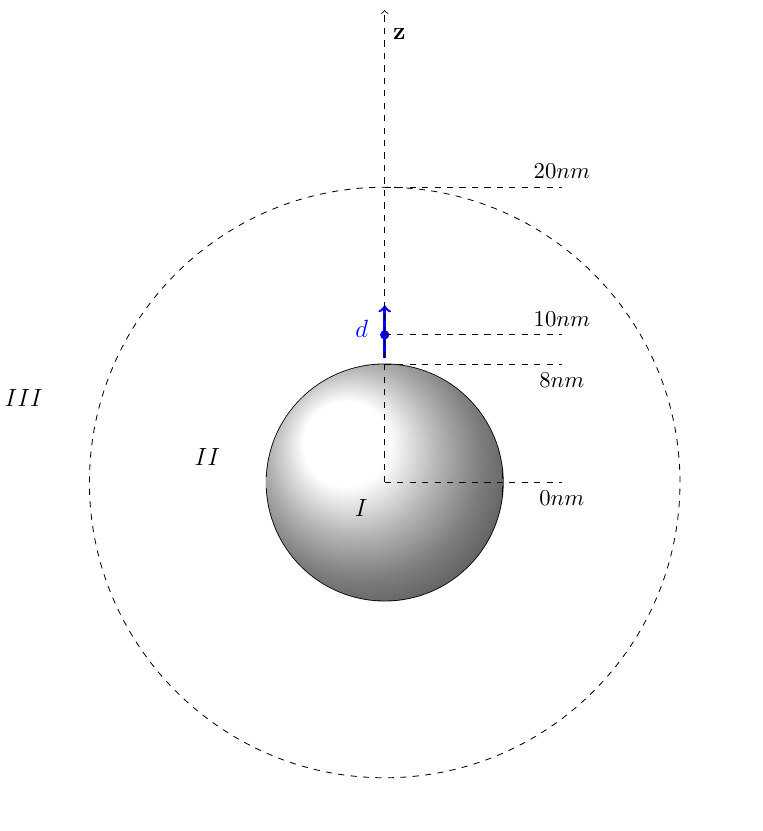}
\caption{On the left: domain boundaries for 3 cigar shaped domains 
all hitting the symmetry axis for $x_{3}=\pm 8$, but the $\varrho$ axis 
in the points 2,4,6 from left to right}; \textcolor{green}{on the right: spherical 
nano-particle of radius 8nm and a dipole on the $x_{3}$ axis.} 
 \label{PS3a}
\end{figure}

For larger values of $\omega$, but the same value of $a$, the resolution has to be adjusted 
since the solution becomes more oscillatory. For 
$\omega=10$, we use $N_{I}=40$, $N_{II}=100$, $N_{III}=30$ and 
$N_{x}=30$ and reach a global error of the order of $10^{-14}$. If 
instead we consider $\omega=1$, but change $a$ and thus the shape of the 
constant coordinate surfaces, no higher resolution is needed. We keep the intersection of the 
cigar like constant coordinate surface such that it intersects the 
$x_{3}$-axis for $x_{3}=8$, but change the 
intersection with the $x_{1},x_{2}$-plane from 6 to 2, this implies 
we use a larger $a\sim 7.75$. The situation is shown in 
Fig.~\ref{PS3a} on the left. With the same numerical parameters as in 
Fig.~\ref{figsolps}, we get again a global error of the order of 
$10^{-14}$. 

%
%

\subsection{Far field of a strongly coupled dipole-nanoparticle system}

In this subsection we study the interaction of a spherical nano-particle 
with a monochromatically radiating dipole. Since we are here 
interested in a concrete physical problem, the speed of light $c$ is 
in this subsection considered in SI units. \\
%
%
The dipole can be placed 
without loss of generality on the $x_{3}$ axis and will be located 
there at $x_{3}=z_{0}$ outside of the nano-particle, but close to its surface. 
It is known that the dipole in vacuum leads to the electric field 
($\mathbf{x_{0}}=z_{0}\hat{x}_{3}$, where $\hat{x}_{3}$ is the unit vector in 
the $x_{3}$ direction)
\begin{equation}
	\begin{split}
	\mathbf{E}_{d}&=\left(\frac{k^{2}}{R^{3}}(\mathbf{(x-x_{0})}\times 
	\mathbf{p})\times 
	\mathbf{(x-x_{0})}\right.\\
	&\left.+(3(\mathbf{x-x_{0}})(\mathbf{(x-x_{0})p})-\mathbf{p}R^{2})\left(\frac{1}{R^{5}}+\frac{ik}{R^{4}}\right)\right)e^{-ikR}
	\end{split}	
	\label{dipole},
\end{equation}
where $k = \omega/c$ and where
\begin{equation}
	R := \sqrt{r^{2}-2rxz_{0}+z_{0}^{2}}.
	\label{Rz0}
\end{equation}
The field $\mathbf{E}_{d}$ is the solution to the vector Helmholtz 
equation (\ref{Helmholtz}) in vacuum,
\begin{equation}
	\nabla\times\nabla\times \mathbf{E}_{d}(\mathbf{x},\omega) - 
k^{2}\mathbf{E}_{d}(\mathbf{x},\omega) =  
\mathbf{S},
	\label{Helmholtzdipole}
\end{equation}
where $\mathbf{S}$ is a distributional source. 
The dipole is chosen to point in the $x_{3}$ direction, 
$\mathbf{p}=p\hat{x}_{3}$. Thus we get 
in spherical coordinates
\begin{equation}
	E^{r}_{d}=pe^{-ikR}\left(\frac{k^{2}z_{0}r(1-x^{2})}{R^{3}}+(3(rx-z_{0})(r-xz_{0})-xR^{2})
	\left(\frac{1}{R^{5}}+\frac{ik}{R^{4}}\right)\right)
	\label{Er}
\end{equation}
and
\begin{equation}
	E^{\theta}_{d}=p\sin\theta e^{-ikR}\left(-\frac{k^{2}r(r-xz_{0})}{R^{3}}+(3(rx-z_{0})z_{0}+R^{2})
	\left(\frac{1}{R^{5}}+\frac{ik}{R^{4}}\right)\right).
	\label{Etheta}
\end{equation}
Since $Y = (\sin\theta (rE^{\theta})_{r}+(1-x^{2})E^{r}_{x})/r$,  we 
get for the twist potential $Y_{d}$ of the dipole
\begin{equation}
	Y_{d} = 
	p(1-x^{2})e^{-ikR}\left(\frac{ik^{3}r}{R^{2}}+\frac{k^{2}r}{R^{3}}\right)
	\label{Ydipole}.
\end{equation}

The radiating dipole to interacts with the spherical metallic
nano-particle. The question is how this nano-particle acts as a 
resonator. We consider a typical situation of strong interaction 
\cite{DRJCCG}, a silver particle of radius of $8nm$, the dipole coupled to it is placed $2nm$ from the north pole on the 
$x_{3}$-axis, thus 
$z_{0}=10nm$, see Fig.~\ref{PS3a} on the right.  The permittivity of the nano-particle is modeled by Drude's formula
\begin{equation}
\epsilon^{I} = \epsilon_{\infty} - 
\frac{\omega_p^2}{\omega^2+\mathrm{i}\gamma_p\omega};  
\end{equation}
for silver \cite{VP12}, the appropriate values are $\epsilon_{\infty} =  6$,
$\hbar\omega_p =  7.90 eV$, and $\hbar\gamma_p = 51 meV$.

Since the dipole is singular at $\mathbf{x}_{0}$ and since 
$\mathbf{S}$ corresponds to a delta function, the problem is not suited for a direct numerical treatment without taking care of the 
singularities. However, we know the dipole field in vacuum, and we 
are only interested in the radiation it causes. Since the Maxwell 
equations and the Helmholtz equation (\ref{Helmholtz}) are 
linear, we can make the ansatz 
$\mathbf{E}=\mathcal{E}+\mathbf{E}_{d}$. With (\ref{Helmholtz}) 
and (\ref{Helmholtzdipole}) we get that the radiation field 
$\mathcal{E}$ satisfies the equation 
\begin{equation}
	\nabla\times\nabla\times \mathcal{E}(\mathbf{x},\omega) - 
\omega^{2}\epsilon(\mathbf{r},\omega)\mathcal{E}(\mathbf{x},\omega) =  	
\omega^{2}(\epsilon(\mathbf{r},\omega)-1)\mathbf{E}_{d}
	\label{HelmholtzE}.
\end{equation}
This means that the dipole field acts as a source for the field 
$\mathcal{E}$, but just where $\epsilon(\mathbf{r},\omega)\neq 1$, in 
our case the interior of the nano-particle. Since the dipole is 
 located outside of the nano-particle, the 
quantity $R$ will not vanish there.  If we introduce the twist potential $\mathcal{Y}$ for 
$\mathcal{E}$, it satisfies with (\ref{NAsp3})
\begin{equation}
	r^{2}(\mathcal{Y}_{rr}+\omega^{2}\epsilon\mathcal{Y})+2r\mathcal{Y}_{r}
	+(1-x^{2})\mathcal{Y}_{xx}=f,
	\label{Yeq}
\end{equation}
where 
\begin{equation}
	f = \omega^{2}(\epsilon -1)
	r\sin\theta((-rE^{\theta}_{d})_{r}+(E^{r}_{d})_{x}) = 
	-\omega^{2}(\epsilon-1) r^{2}Y_{d}.
	\label{fdipole}
\end{equation}
Thus
the source $f$  vanishes for $r\to0$. We normalize the solution in 
the following at infinity such that in the $x_{3}=0$ plane we have
\begin{equation}
	\lim_{r\to\infty}\frac{\mathcal{Y}}{Y_{d}}=1.
	\label{Ynorm}
\end{equation}
In other words we divide $\mathcal{Y}$ by $ik^{3}p$. 

\begin{figure}[htb!]
  \includegraphics[width=0.32\textwidth]{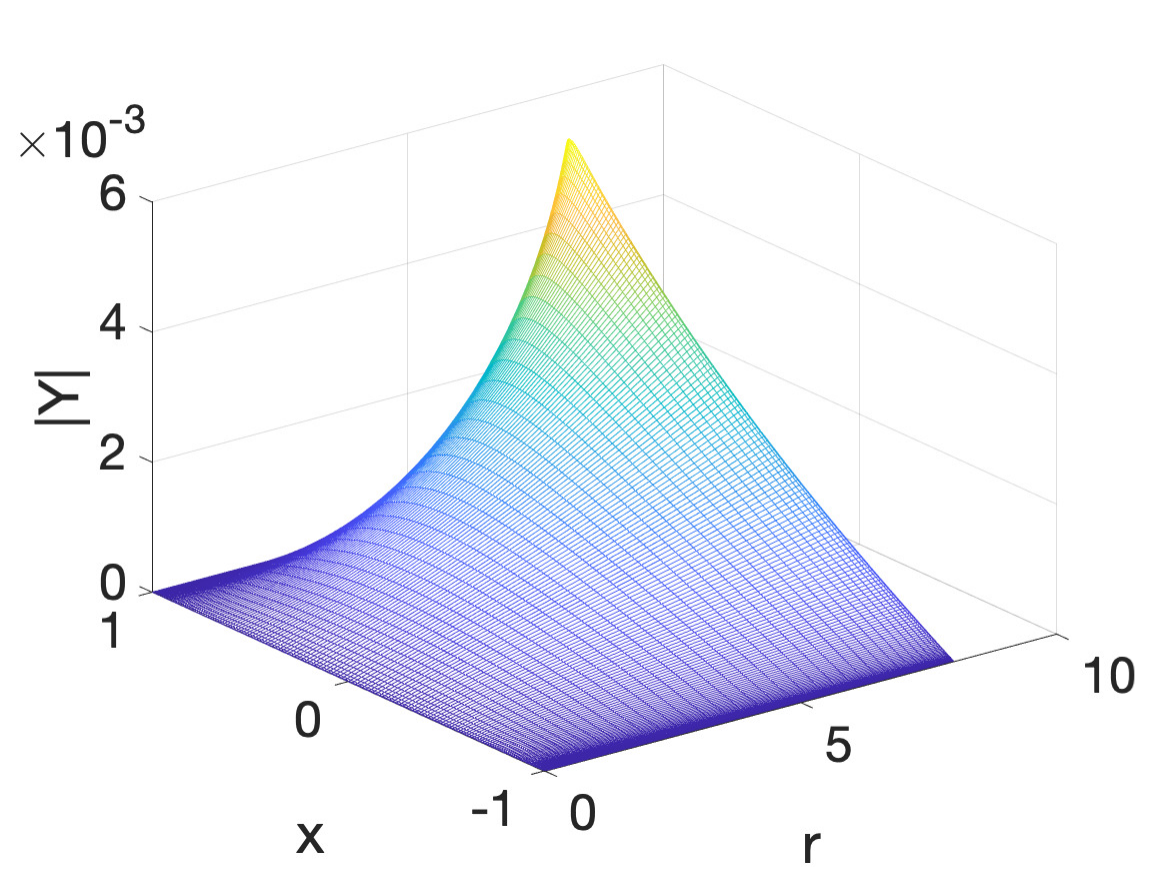}
  \includegraphics[width=0.32\textwidth]{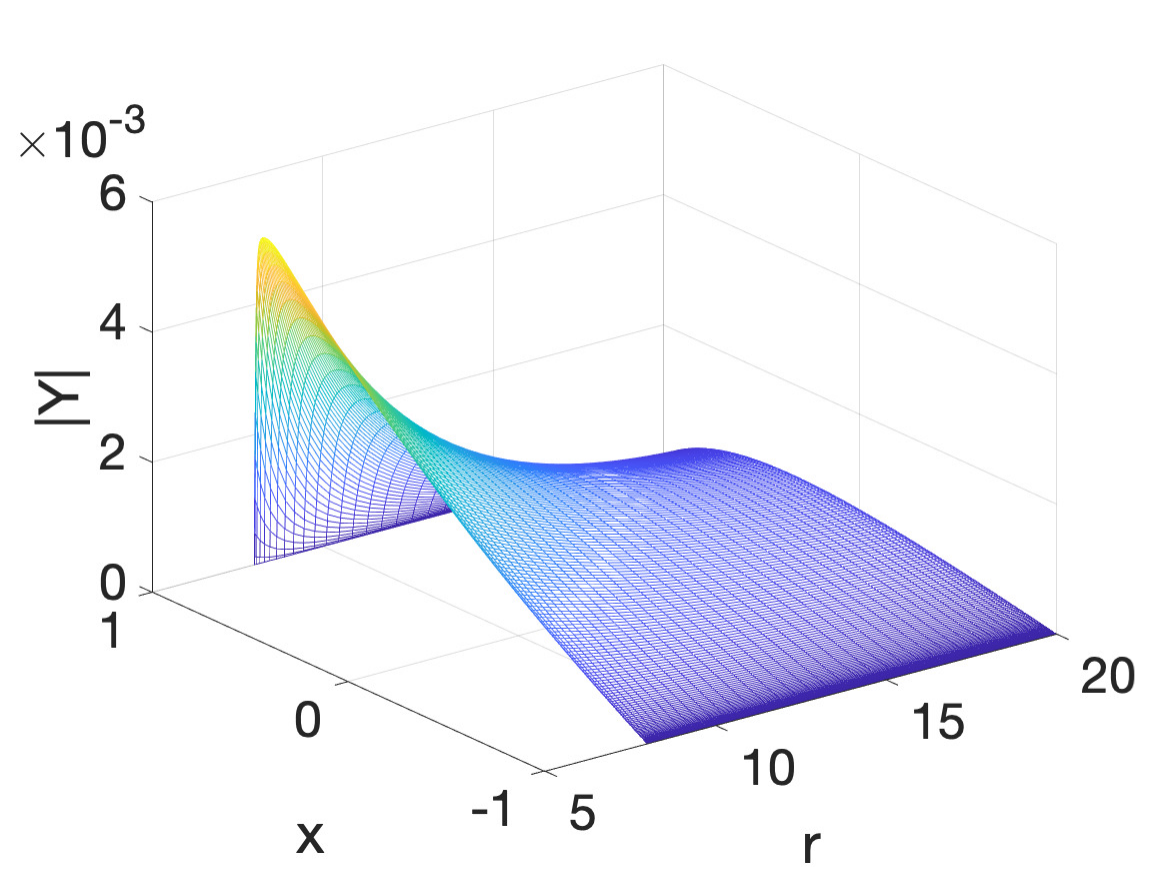}
  \includegraphics[width=0.32\textwidth]{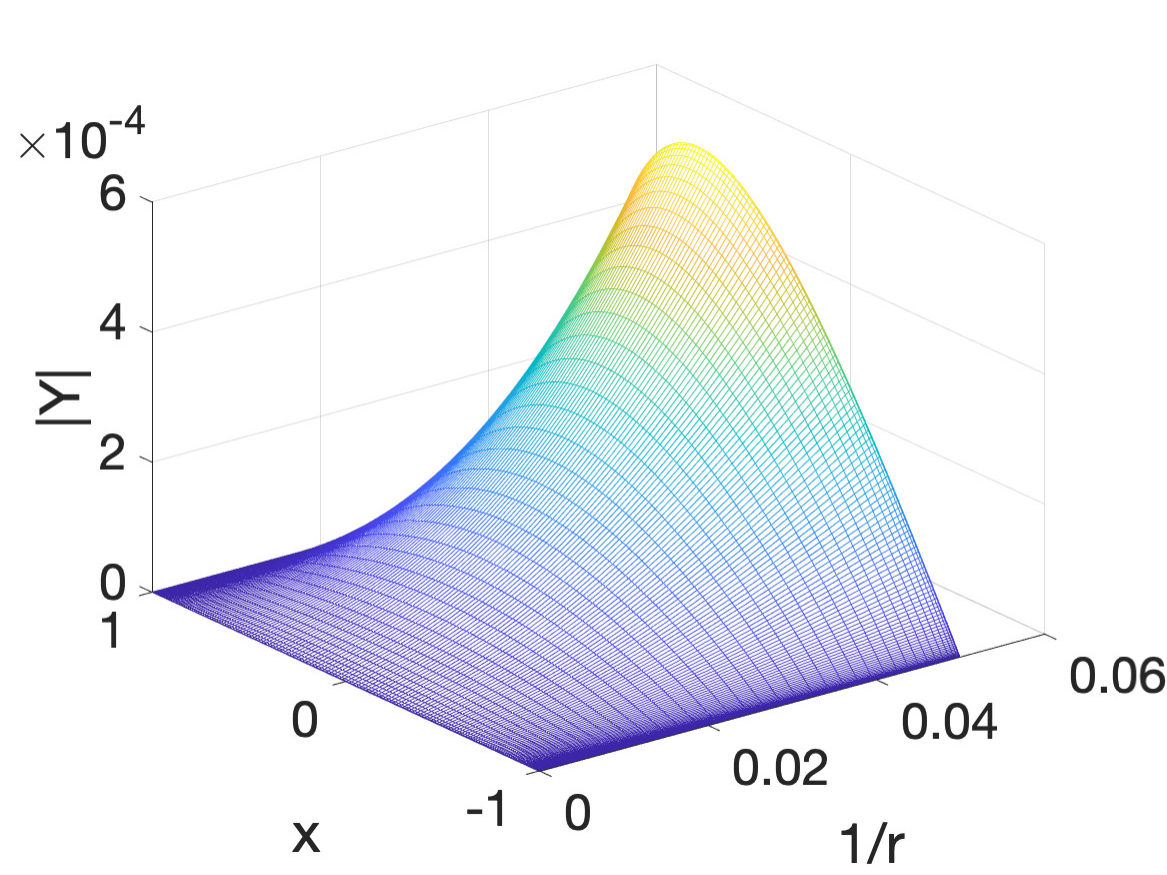}\\
  \includegraphics[width=0.32\textwidth]{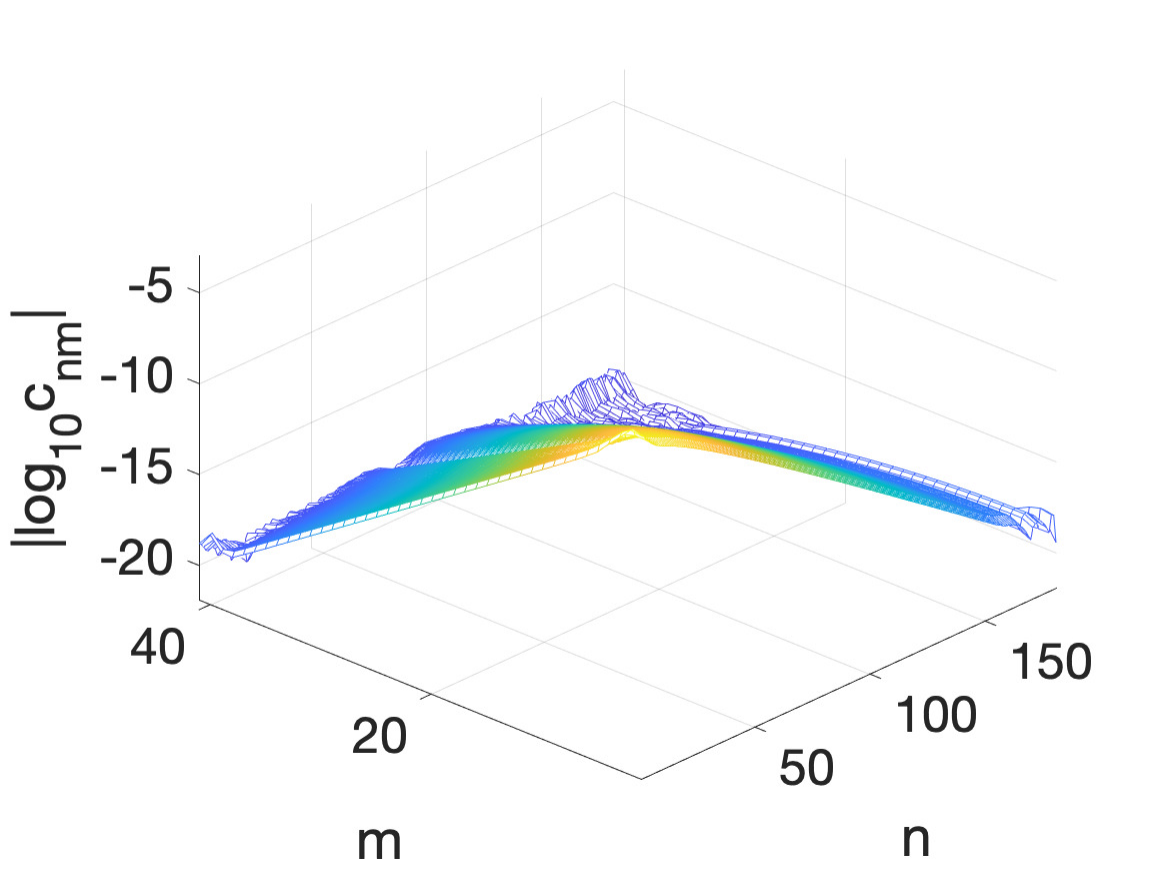}
  \includegraphics[width=0.32\textwidth]{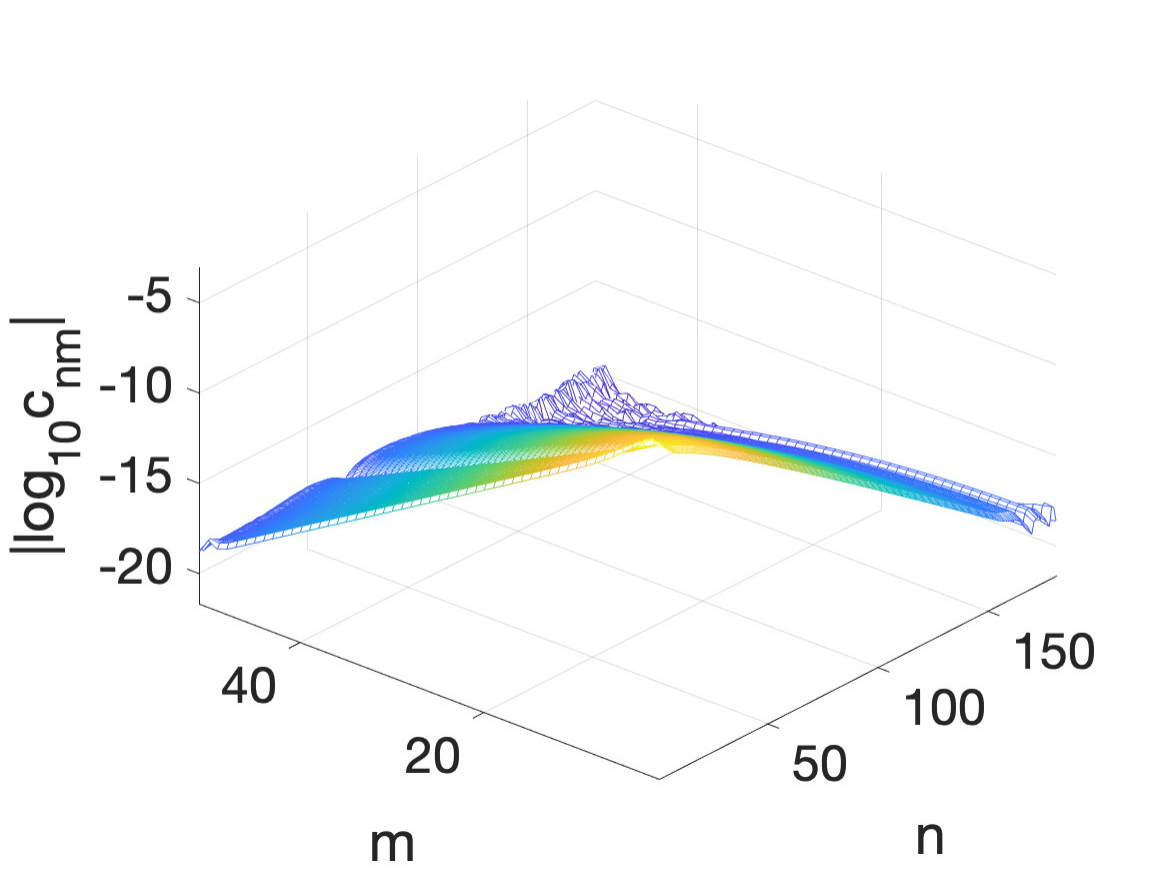}
  \includegraphics[width=0.32\textwidth]{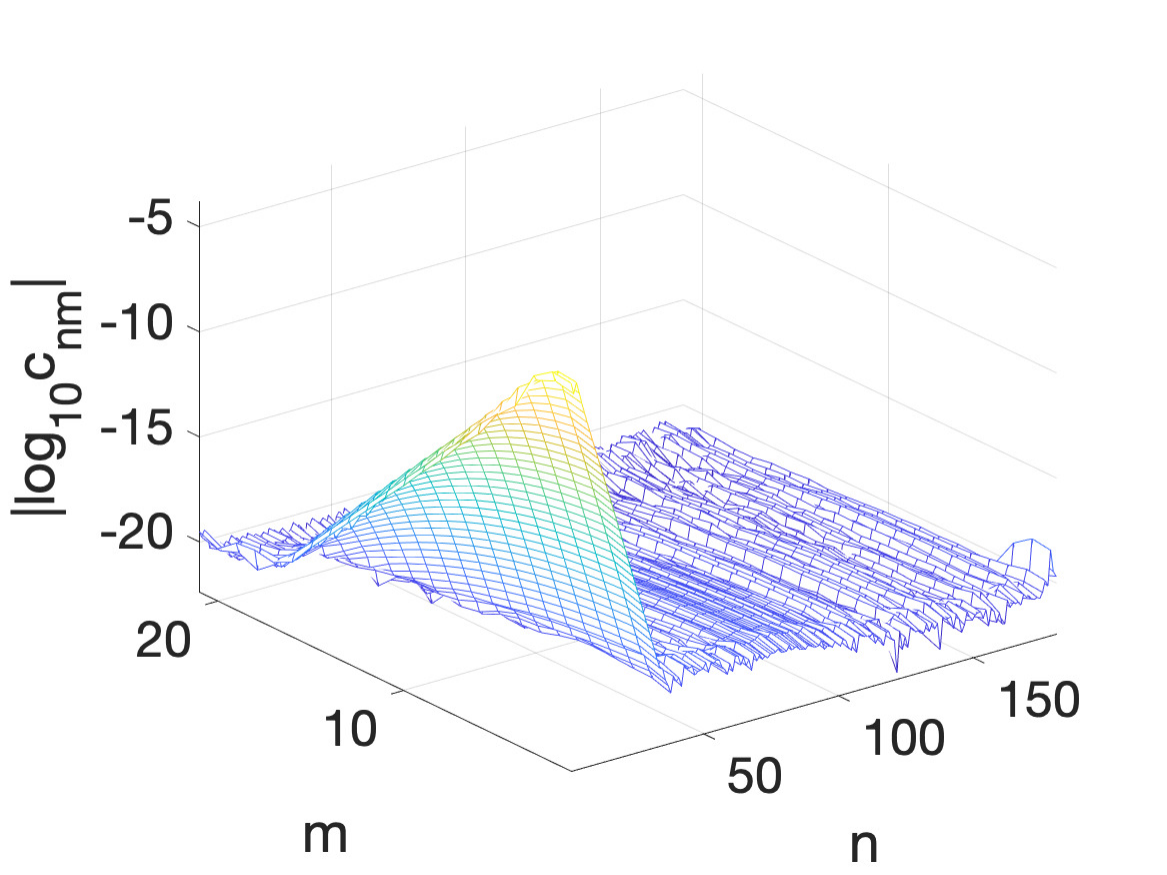}
\caption{Solution to the equation (\ref{Yeq}) for $\omega=3eV/\hbar$ in the 
 domains I, II, III (from left to right) in the upper row, and the 
 corresponding
 spectral coefficients in the lower row.}
 \label{Yom3dipole}
\end{figure}

The resolution for the computation of $\mathcal{Y}$ can be 
estimated from the spectral coefficients of $f$ (they have to decrease to the order of machine 
precision, for double precision this is roughly $10^{-16}$) in domain I. We work in the following with the domains as 
before, I ($r < 8$), II ($8 < r < 20$) and III ($r > 20$) and use $N_{I}=40$, 
$N_{II}=50$ and $N_{III}=20$ and $N_{x}=180$. 
The solution for 
$\omega=3eV/\hbar$ (a typical frequency considered in this context in 
nano-optics) can be seen in Fig.~\ref{Yom3dipole}. The spectral 
coefficients of the solution in each domain in the lower row of the 
figure indicate that spectral accuracy is achieved.

\begin{figure}[htb!]
  \includegraphics[width=0.7\textwidth]{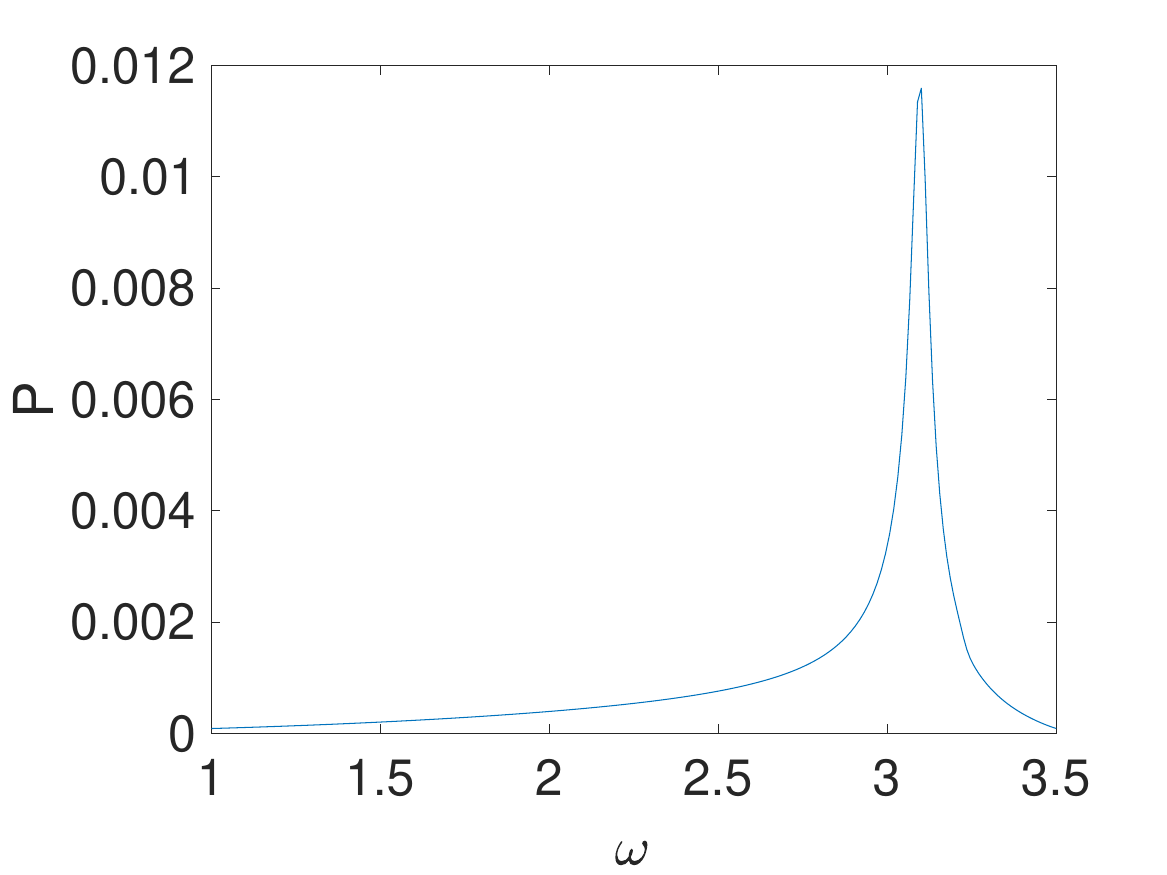}
 \caption{Induced dipole field in the nano-particle normalized by the 
 dipole field at infinity in dependence of $\omega$.}
 \label{Pdipole}
\end{figure}

For a dipole the strongest radiation is always expected orthogonal to 
the orientation of the dipole, here in the far field in the $\{x_{1}, x_2\}$ plane. The induced 
dipole moment in the nanoparticle can be read of at infinity. Since 
the leading contribution will be linear in $1/r$ in this case, we 
simply differentiate $\mathcal{Y}$ in domain III with respect to 
$\rho=1/r$ and denote by $P$ the resulting value divided by $ik^{3}p$, 
the leading contribution of the dipole (\ref{Ydipole}) at infinity. In Fig.~\ref{Pdipole} we plot this quantity for several 
values of $\omega$.  The strongest resonance of the nano-particle is 
observed for $\omega\approx 3eV/\hbar$.


\section{Conclusion}
In this paper we have presented a multi-domain spectral approach for 
the monochromatic Maxwell equations in an 
axisymmetric setting in spherical and prolate spheroidal 
coordinates. The Sommerfeld condition is imposed as in 
\cite{BM81,GJ17,GTZ18} exactly at infinity after splitting off an 
oscillatory factor. For several examples it is shown that machine 
precision can be reached with this approach. Obviously one could have 
built a similar spectral approach based on the eigenfunctions of the 
Helmholtz equation, in spherical coordinates Legendre polynomials and 
spherical Bessel functions. The differentiation matrices for Legendre 
polynomials are known, see for instance \cite{trefethen}, but the 
spherical Bessel functions are transcendental functions that have to 
be
computed as well. In particular the treatment of the Hankel function 
at infinity will need a similar treatment as presented here. In 
contrast to the case of Chebyshev polynomials, no fast algorithm to 
compute the spectral coefficients is known. The situation is worse in 
the prolate spheroidal case where the eigenfunctions are less well 
known. Thus it appears that the numerical method we discuss here 
could be also suitable to efficiently compute these functions which 
will be studied elsewhere along the lines of \cite{hypergeom}.

Whereas we consider in this paper only the axisymmetric case, the 
approach is set up in a way that it can be extended to situations 
without symmetry. Instead of one twist potential, one has to 
deal with all three components of the electric field in this case in 
the Helmholtz equation (\ref{Helmholtz}). The discretization in $r$ 
and $x$ is as presented above. The dependence on the 
azimuthal coordinate $\phi$ can be addressed with a Fourier spectral 
method (see for instance \cite{trefethen}) which has the advantage of 
diagonal differentiation matrices. This means the equations decouple 
in $\phi$.
For each of the $N_{\phi}$ collocation points in $\phi$, one thus has 
to solve a system for $\mathbf{E}$ with the methods discussed in the 
present 
paper. Since the equations do not couple in $\phi$, this is fully 
parallizable.  A full 3D code along these lines will be the subject of further work.

\end{document}